\documentclass[preprint,11pt]{elsarticle}
\usepackage{color,amsmath,amssymb,graphicx,algorithmic,alltt,xy,fullpage,amsthm,wrapfig,setspace,cancel}

\definecolor{refkey}{rgb}{1,0.5,0.5}
\definecolor{labelkey}{rgb}{0.5,1,0.5}
\definecolor{MyDarkGreen}{rgb}{0,0.45,0}
\definecolor{MyDarkBlue}{rgb}{0,0,0.75}
\definecolor{MyDarkRed}{rgb}{0.9,0,0}
\date{}
\input xy
\xyoption{all}
\newtheorem{thrm}{Theorem}
\newtheorem{lem}{Lemma}

\newtheorem{rmk}{Remark}

\newtheorem{conc}{Conclusion}
\newtheorem{corol}{Corollary}
\newtheorem{prop}{Continuation Result}

\newcommand{\X} {{\bf x }}
\newcommand{\x} {{\bf x }}
\newcommand{\schro} {Schr\"odinger}

\newcommand{\Beginproof}{\vspace{0mm} \parindent=0pt
         {\bf Proof.} \hspace{3mm} \parindent=3ex}
\newcommand{\Endproof}{$\Box$ \vspace{5mm}
                        \parindent=3ex}
\newcommand{\psiex}{\psi_{\rm explicit}}
\newcommand{\psiexalpha}{\psi_{\rm explicit, \alpha}}
\newcommand{\Pcr}{P_{\rm cr}}

\newcommand{\rev} { }

\DeclareMathOperator*{\argmax}{arg\,max}

\graphicspath{{}}
\begin{document}
\begin{frontmatter}
\title{Continuations of the nonlinear Schr\"odinger equation beyond the singularity}

\author{G. Fibich}\ead{fibich@tau.ac.il}\author{M. Klein}\ead{morankli@tau.ac.il}
\address{School of Mathematical Sciences, Tel Aviv University, Tel Aviv 69978, Israel}
\begin{abstract}
We present four continuations of the critical nonlinear
\schro~equation (NLS) beyond the singularity: 1) a sub-threshold
power continuation, 2) a shrinking-hole continuation for ring-type
solutions, 3) a vanishing nonlinear-damping continuation, and 4) a
complex Ginzburg-Landau (CGL) continuation. Using asymptotic
analysis, we explicitly calculate the limiting solutions beyond the
singularity. These calculations show that for generic initial data
that leads to a loglog collapse, the sub-threshold power limit is a
Bourgain-Wang solution, both before and after the singularity, and
the vanishing nonlinear-damping and CGL limits are a loglog solution
before the singularity, and have an infinite-velocity{\rev{expanding
core}} after the singularity.  Our results suggest that all NLS
continuations share the universal feature that after the singularity
time~$T_c$, the phase of the singular core is only determined up to
multiplication by~$e^{i\theta}$. As a result, interactions between
post-collapse beams (filaments) become chaotic. We also show that
when the continuation model leads to a point singularity and
preserves the NLS invariance under the
transformation~$t\rightarrow-t$ and~$\psi\rightarrow\psi^\ast$, the
singular core of the weak solution is symmetric with respect
to~$T_c$. Therefore, the sub-threshold power and
the{\rev{shrinking}}-hole continuations are symmetric with respect
to~$T_c$, but continuations which are based on perturbations of the
NLS equation are generically asymmetric.

\end{abstract}
\end{frontmatter}

\section{Introduction}\label{sec:Introduction} The
focusing nonlinear Schr\"odinger equation (NLS)
\begin{equation}\label{eq:DdimensionalNLS}
i\psi_t(t,\X)+\Delta\psi+|\psi|^{2\sigma}\psi=0,\qquad\psi_0(0,\X)=\psi_0(\X)\in
H^1,
\end{equation}
where~$\X=(x_1,...,x_d)\in\mathbb{R}^d$
and~$\Delta=\partial_{x_1x_1}+\cdot\cdot\cdot\partial_{x_dx_d}$, is
one of the canonical nonlinear equations in physics, arising in
various fields such as nonlinear optics, plasma physics,
Bose-Einstein condensates (BEC), and surface waves. In the
two-dimensional cubic case, this equation models the propagation of
intense laser beams in a bulk Kerr medium. In that case, $\psi$ is
the electric field envelope, $t$~is the direction of propagation,
$d=2$, $x_1$ and~$x_2$ are the transverse coordinates, and
$\sigma=1$ (cubic nonlinearity).

In 1965, Kelley showed that the two-dimensional cubic NLS
admits solutions that collapse (become singular) at a finite
time (distance)~$T_c$~\cite{Kelley-65}. Since physical quantities do not become
singular, this implies that some of the terms that were neglected in
the derivation of the NLS, become important near the singularity.
Therefore, the standard approach for continuing the solution beyond
the singularity has been to consider a more comprehensive model, in
which the collapse is arrested.

In this study, we adopt a different approach, and ask whether
singular NLS solutions can be continued beyond the singularity, {\em
within} the NLS model. By this we mean that the solution satisfies
the NLS both before and after the singularity, and a matching
("jump") condition at the singularity. The motivation for this
approach comes from hyperbolic conservation laws, where in the
absence of viscosity, the solution can become singular (develop
shocks). In that case, there is a huge body of literature on how to
continue the inviscid solution beyond the singularity, which
consists of Riemann problems, vanishing-viscosity solutions, entropy
conditions, Rankine-Hogoniot jump conditions, specialized numerical
methods, etc. In contrast, two studies from 1992 by
Merle~\cite{Merle-92a,Merle-92b}, and a recent study by Merle,
Raphael and Szeftel~\cite{Merle-Raphael-Szeftel}, addressed this
question in the NLS.{\rev{Tao~\cite{Tao-2009} proved the global
existence and uniqueness in the semi Strichartz class for solutions
of the critical NLS. Intuitively, these  solutions are formed by
solving the equation in the Strichartz class whenever possible, and
deleting any power that escapes to spatial or frequency infinity
when the solution leaves the Strichartz class. These solutions,
however, do not depend continuously  on the initial conditions, and
are thus not a well-posed class of solutions.
Stinis~\cite{Stinis-2010} studied numerically the continuation of
singular NLS solutions using the t-model approach.}}

 In~\cite{Merle-92a}, Merle presented an explicit continuation of a singular NLS solution
beyond the singularity, which is based on reducing the power ($L^2$
norm) of the initial condition of the explicit blowup
solution~$\psiexalpha$,{\rev{see~(\ref{eq:ExplicitBlowupSolutionAlpha})}},
below the critical power for collapse~$\Pcr$. This continuation has
two key properties:
\begin{enumerate}
\item \textbf{Property 1}: The solution is symmetric with respect to the singularity
time~$T_c$.\label{Item:WeakSolProp1}
\item \textbf{Property 2}: After the singularity, the solution can only be determined up to multiplication by a constant
phase term.\label{Item:WeakSolProp2}
\end{enumerate}
Merle's breakthrough result, however, applies only to the critical
NLS ($\sigma d=2$), and only to the explicit blowup
solutions~$\psiexalpha$.
Recently, Merle, Raphael and Szeftel~\cite{Merle-Raphael-Szeftel} generalized this result to
Bourgain-Wang singular solutions~\cite{Bourgain-97}, i.e., solutions that have a
singular component that collapses as~$\psiexalpha$, and a
non-zero regular component that vanishes at the singularity
point and does not participate in the collapse process.

In~\cite{Merle-92b}, Merle presented a different continuation, which
is based on the addition of nonlinear saturation. This study showed
that, generically, as the nonlinear saturation parameter goes to
zero, the limiting solution can be decomposed beyond~$T_c$ into two
components: A $\delta$-function singular component with power~$m(t)
\ge \Pcr$, and a regular component elsewhere. Similar results follow
from the asymptotic analysis of Malkin~\cite{Malkin-1993}, which
suggests that $m(t) \equiv  \Pcr$ for $T_c \le t <\infty$.

It is thus useful to distinguish between two types of continuations:
\begin{enumerate}
  \item A {\bf point singularity}, in which the weak solution is singular at~$T_c$, but regular (i.e., in~$H^1$) for $t>T_c$.
  \item A {\bf filament singularity}, in which the weak solution has a $\delta$-function singularity for $T_c \le t \le T_0$, where $T_c<T_0 \le \infty$.
  The motivation for this terminology comes from nonlinear optics,
  where collapsing beams can form long and narrow filaments (which,
  to leading order, can be viewed as an "extended"~$\delta$-function).
\end{enumerate}

In this work we propose four novel continuations that lead to a
point singularity, and obtain explicit formulae for the solution
beyond the singularity. Our main findings are:
\begin{enumerate}
\item The non-uniqueness of the phase beyond the singularity (Property~2) is a
universal feature of{\rev{NLS}} continuations.
\item The symmetry with respect to the singularity
time (Property~1) holds only if the continuation is time reversible and leads to a point singularity.
Therefore, it is non-generic.
\end{enumerate}

The paper is organized as follows. In
section~\ref{sec:NLS_TheoryReview} we provide a short review of NLS
theory. In section~\ref{sec:WeakSolNLS} we present Merle's
continuation of~$\psiexalpha$, and illustrate it numerically. In
section~\ref{sec:NonExplicitBeam} we generalize this approach and
present a sub-threshold power continuation, which can be applied
to{\rev{generic}} initial condition of the critical NLS. We compute
asymptotically the limiting solution, and show that it is a
Bourgain-Wang solution, both before and after the singularity. In
particular, this continuation preserves the two key properties of
Merle's continuation. In section~\ref{sec:SymmetryProperty} we show
that Property~1 (symmetry with respect to~$T_c$) holds for
time-reversible continuations that lead to a point singularity. In
Section~\ref{sec:nl-saturation} we discuss the nonlinear-saturation
continuation, which leads to a filament singularity. In
section~\ref{sec:InteractionTwoBeams} we show that because of the
phase non-uniqueness (Property~2), the
interaction between post-collapse beams is chaotic. In
section~\ref{sec:RingTypeSingularSolution} we present a
vanishing-hole continuation, which is suitable for ring-type
singular solutions. This continuation is time-reversible, and it
satisfies Properties~1 and~2. In section~\ref{sec:DampedNLS} we
present a vanishing nonlinear-damping continuation, and compute the
continuation asymptotically in two cases:
\begin{enumerate}
  \item The  nonlinear-damping continuation of the explicit solution~$\psiexalpha$
is, up to an undetermined constant phase, given by~$\psi_{\rm explicit, \kappa \alpha}$ with $\kappa \approx 1.614$.

  \item {\rev{The nonlinear-damping continuation of solutions that undergo a loglog collapse has an infinite-velocity expanding core,
  with an undetermined constant phase.}}
\end{enumerate}
 Therefore, the phase becomes non-unique beyond the
singularity (Property~2). In contrast with
previous continuations, however, the solution is asymmetric with
respect to the singularity time (i.e., the solution does not satisfy
Property~1). This is to be expected, as the
nonlinear-damping continuation is not time-reversible. In
Section~\ref{sec:CGL} we show that the continuation which is based
on the complex Ginzburg-Landau (CGL) limit of the NLS, is equivalent
to the vanishing nonlinear-damping continuation. In
section~\ref{sec:SingularityInLS} we show a continuation of singular
solutions of the linear \schro~equation. In this case, the limiting
phase beyond the singularity is unique. This shows that the
post-collapse non-uniqueness of the phase
(Property~2) is a nonlinear phenomenon.
Section~\ref{sec:FinalRemarks} concludes with a discussion.

\subsection{Level of rigor}

The results which are derived in this manuscript are non-rigorous, and are based on
asymptotic analysis, numerical simulations, and physical arguments.
To emphasize this, we use the terminology {\em Continuation Results}, rather than Propositions or Theorems.

\section{Review of NLS theory}\label{sec:NLS_TheoryReview}

We briefly review NLS theory, for more information
see~\cite{PNLS-99, Sulem-99, Strauss-89}. The
NLS~(\ref{eq:DdimensionalNLS}) has two important conservation laws:
{\em {Power}} conservation\footnote{We call the~$L^2$ norm the
{\em{power}}, since in optics it corresponds to the beam's power.}
\[
P(t)\equiv P(0),\qquad P(t)=\int|\psi|^2d\X,
\]
and {\em{Hamiltonian}} conservation
\begin{equation}\label{eq:HamiltonianConservation}
H(t)\equiv H(0),\qquad
H(t)=\int|\nabla\psi|^2d\X-\frac1{\sigma+1}\int|\psi|^{2\sigma+2}d\X.
\end{equation}

The NLS admits the waveguide solutions~$\psi=e^{it}R(r)$,
where~$r=|\x|$, and~$R$ is the solution of
\begin{equation}\label{eq:dDim_R_ODE}
R''(r)+\frac{d-1}{r}R'-R+R^{2\sigma+1}=0,\qquad R'(0)=0,\quad
R(\infty)=0.
\end{equation}
When~$d=1$, the solution of~(\ref{eq:dDim_R_ODE}) is unique, and is
given by
\begin{equation}\label{eq:1D_ground_state_solutions}
R_\sigma(x)=(1+\sigma)^{1/2\sigma}\mbox{cosh}^{-1/\sigma}(\sigma x).
\end{equation}
When~$d\geq2$, equation~(\ref{eq:dDim_R_ODE}) admits an infinite
number of solutions. The solution with the minimal power, which we
denote by~$R^{(0)}$, is unique, and is called the ground state.

When~$\sigma d<2$, the NLS is called subcritical. In that case,
all~$H^1$ solutions exist globally. In contrast, both the critical
NLS ($\sigma d=2$) and the supercritical NLS ($\sigma d>2)$ admit
singular solutions.

Let~$\psi(t,\bf x)$ be a solution of the
NLS~(\ref{eq:DdimensionalNLS}). Then,~$\psi$ remains a solution of
the NLS~(\ref{eq:DdimensionalNLS}) under the following
transformations:
\begin{enumerate}
\item Spatial translations:~$\psi(t,\x)\rightarrow\psi(t,\x+\x_0)$,
where~$\x_0\in\mathbb{R}^d$. \label{item:Spatial translations}
\item Temporal translations:~$\psi(t,\X)\rightarrow\psi(t+t_0,\X)$,
where~$t_0\in\mathbb{R}$.\label{item:Axial translation}
\item Phase change:~$\psi(t,\x)\rightarrow e^{i\theta}\psi(t,\x)$,
where~$\theta\in\mathbb{R}$. \label{item:Phase change}
\item
Dilation:~$\psi(t,\X)\rightarrow\lambda^{1/\sigma}\psi(\lambda^2t,\lambda\X)$,
where~$\lambda\in\mathbb{R}^+$.\label{item:Dilation}
\item Galilean
transformation:~$\psi(t,\x)\rightarrow\psi(t,\x-{\bf c} t)e^{i{\bf
c}\cdot\x/2-i|{\bf c}|^2t/4}$, where~${\bf c}\in\mathbb{R}^d$.
\label{item:Galilean transformation}
\end{enumerate}
Therefore, multiplying the initial condition by a constant
phase~$e^{i\theta}$ does not affect the solution. In addition, by
the Galilean transformation, multiplying the initial condition by a
linear phase term~$\psi_0(x)\rightarrow\psi_0(\X)e^{i{\bf
c}\cdot\X/2}$ does not affect the dynamics, but rather causes the
solution to be tilted in the direction of~${\bf
{n}}=(1,{\bf{c}})\in\mathbb{R}\times\mathbb{R}^d$.
\subsection{Critical NLS}
In the critical case~$\sigma d=2$,
equation~(\ref{eq:DdimensionalNLS}) can be rewritten as
\begin{equation}
   \label{eq:dDim_CriticalNLS}
i\psi_t(t,\X)+\Delta\psi+|\psi|^{4/d}\psi=0,\qquad\psi_0(0,\X)=\psi_0(\X)\in
H^1,
\end{equation}
and equation~(\ref{eq:dDim_R_ODE}) can be rewritten as
\begin{equation}\label{eq:dDimCritical_R_ODE}
R''(r)+\frac{d-1}{r}R'-R+R^{4/d+1}=0,\qquad R'(0)=0,\quad
R(\infty)=0.
\end{equation}
\begin{thrm}[Weinstein~\cite{Weinstein-83}] \label{Theorem:WeinsteinThm}
 A sufficient condition for global
existence in the critical NLS~(\ref{eq:dDim_CriticalNLS}) is
that~$\|\psi_0\|_2^2<\Pcr$, where~$\Pcr=\|R^{(0)}\|_2^2$,
and~$R^{(0)}$ is the ground state of
equation~(\ref{eq:dDimCritical_R_ODE}).
\end{thrm}

The critical NLS~(\ref{eq:dDim_CriticalNLS}) admits the explicit
solution
\begin{subequations}
 \label{eq:ExplicitBlowupSolution}
\begin{equation}\label{eq:ExplicitBlowUpSolutions_1}
\psiex(t,r)=\frac1{L^{d/2}(t)}R^{(0)}\left(\frac{r}{L(t)}\right)e^{i\tau+i\frac{L_t}{L}\frac{r^2}{4}},
\end{equation}
where
\begin{equation}\label{eq:ExplicitBlowupSolution_2}
L(t)=T_c-t,\qquad \tau(t)=\int_0^t\frac1{L^2(s)} \, ds=\frac1{T_c-t}.
\end{equation}
\end{subequations}
More generally, applying the dilation transformation with $\lambda = \alpha$
and the temporal translation $T_c\longrightarrow \alpha^2 T_c$ shows that
the critical NLS~(\ref{eq:dDim_CriticalNLS}) admits
the explicit solutions
\begin{subequations}\label{eq:ExplicitBlowupSolutionAlpha}
\begin{equation}\label{eq:ExplicitBlowupSolutionAlpha_1}
\psiexalpha(t,r)=\frac1{L_\alpha^{d/2}(t)}R^{(0)}\left(\frac{r}{L_\alpha(t)}\right)e^{i\tau_\alpha+i\frac{(L_\alpha)_t}{L_\alpha}\frac{r^2}{4}},
\end{equation}
where
\begin{equation}\label{eq:ExplicitBlowupSolutionAlpha_2}
L_\alpha(t)=\alpha(T_c-t),\qquad
\tau_\alpha(t)=\int_0^t\frac1{L_\alpha^2(s)} \, ds=\frac1{\alpha^2}\frac1{T_c-t},\qquad\alpha>0.
\end{equation}
\end{subequations}
Even more generally, by the Galilean transformation,
\begin{equation}\label{eq:GeneralExplicitBlowupSol}
\psiexalpha^{tilt,{\bf{c}}}(t,\X)=\psiexalpha(t,\X-{\bf
c}\cdot t)e^{ i{\bf c}\cdot\X/2-i|{\bf c}|^2t/4},
\end{equation}
are also explicit solutions of the critical
NLS~(\ref{eq:dDim_CriticalNLS}).

The explicit
solutions~(\ref{eq:ExplicitBlowupSolution})--(\ref{eq:GeneralExplicitBlowupSol})
become singular at~$t=T_c$. These solutions are unstable, however,
as they have exactly the critical power for collapse. Therefore, any
infinitesimal perturbation which decreases their power, will arrest
the collapse.

The only minimal-power blowup solutions (i.e., singular solutions
whose power is exactly~$\Pcr$) are given
by~$\psiexalpha^{tilt,{\bf{c}}}$:
\begin{thrm}[Merle~\cite{Merle-92a,Merle-92b}]
   \label{Theorem:minimal_power}
Let $\psi$ be a solution of the critical
NLS~(\ref{eq:dDim_CriticalNLS}) which blows up at a finite time
$T_c>0$, such that $||\psi_0||_2^2 =  \Pcr$. Then, there exist
$\alpha\in\mathbb{R}^+$,~$\theta \in\mathbb{R}$, and ${\bf x}_0,
{\bf c} \in \mathbb{R}^d$, such that for $0 \le t<T_c$,
\begin{equation}
   \label{eq:minimal_power_generic_general}
\psi(t, {\bf x}) =
\psiexalpha^{tilt,{\bf{c}}}(t,\X-\X_0)e^{i\theta}.
\end{equation}
\end{thrm}

When an NLS solution whose power is slightly above~$\Pcr$ undergoes a stable
collapse, it splits into two components: A collapsing core that
approaches the universal $\psi_{R^{(0)}}$~profile and blows up at
the loglog law rate, and a non-collapsing tail~($\phi$) that does
not participate in the collapse process:
\begin{thrm}[Merle and Raphael~\cite{Merle-Raphael-2003}, \cite{Merle-04}, \cite{Merle-05}, \cite{Merle-05b},
\cite{Merle-06}, \cite{Merle-06b},\cite{Raphael-05}]
\label{thm:Merle+Raphael}
 Let $d=1,2,3,4,5$, and let~$\psi$ be a
solution of the critical NLS~(\ref{eq:dDim_CriticalNLS}) that
becomes singular at~$T_c$. Then, there exists a universal
constant~$m^*>0$, which depends only on the dimension, such that for
any $\psi_0 \in H^1$ such that
$$
\Pcr <\|\psi_0\|_2^2< \Pcr+m^*, \qquad
   H_G(\psi_0):= H(\psi_0)-\left ( \frac{\mbox{Im} \int \psi_0^* \nabla \psi_0  }{||\psi_0||_2} \right)^2<0,
$$
the following hold:
\begin{enumerate}
  \item    {\rev{There exist parameters
  $(\tau(t), {\bf x}_0(t), L(t)) \in {\Bbb R} \times{\Bbb R}^d\times {\Bbb R}^+$,
  and a function $0\neq\phi \in L^2$,  such that
$$
 \psi(t,{\bf x}) -  \psi_{R^{(0)}}(t,{\bf x}-  {\bf x}_0(t))
 \stackrel{L^2}{\longrightarrow} \phi ({\bf x}), \qquad t \longrightarrow T_c,
$$
where
$$
\psi_{R^{(0)}}(t,{\bf x})  = \frac{1}{L^{d/2}(t)}  R^{(0)} \left(\frac{|{\bf x}|}{L(t)} \right)e^{i \tau(t)},
$$
and $R^{(0)}$ is the ground state of
equation~(\ref{eq:dDimCritical_R_ODE}).}}
\item  As $t \longrightarrow T_c$,
\begin{equation}
\label{eq:loglog_MR} L(t)  \sim \sqrt{2 \pi}
\left(\frac{T_c-t}{\log|\log(T_c-t)|}\right)^{1/2} \!\!\!\!\!\!\!\!
\qquad\mbox{\bf (loglog law)}.
\end{equation}
\end{enumerate}
\end{thrm}

NLS solutions whose power is slightly above~$\Pcr$ can also undergo
a different type of collapse, in which the collapsing core
approaches the explicit blowup solution~$\psiexalpha$ and blows up
at a linear rate,{\rev{but}} the solution also has a nontrivial
tail~($\phi$) that does not participate in the collapse process:
\begin{thrm}[Bourgain and Wang~\cite{Bourgain-97}] Let $d=2$, let~$A_0$ be a given integer, and let~$A\geq A_0$ be a
large enough integer. Let $\phi\in X_A = \{f\in H^A~with~(1+|x|^A)f\in L^2\}$,
and let~$z\in\mathcal{C}((T^\ast,T_c],X^A)$ be the solution to the critical NLS~(\ref{eq:dDim_CriticalNLS}),
subject to $z(t=T_c)=\phi$,
where~$T^\ast<T_c$ is the maximal time of existence of~$z$. Assume
that~$\phi$ vanishes to high order at the origin, i.e.,
$D^\alpha \phi(0)=0$ for $|\alpha|\leq A-1$.
Then,  
there exists~${\rev{T^\ast<}}t_0<T_c$ and a unique
solution~$\psi_{\rm BW}\in\mathcal{C}([t_0, T_c),X_{A_0})$
to~(\ref{eq:dDim_CriticalNLS}), such that
\begin{equation}\label{eq:BourgainWang_3}
\|\psi_{\rm BW}(t)-\psiexalpha(t)-z(t)\|_{X_{A_0}}\leq|T_c-t|^{A_0}.
\end{equation}
\end{thrm}

The Bourgain-Wang solutions~$\psi_{\rm BW}(t,{\bf x})$ are unstable~\cite{Merle-Raphael-Szeftel},
since their singular core is the unstable blowup solution~$\psiexalpha$.

\section{Merle's first continuation}
 \label{sec:WeakSolNLS}

In \cite{Merle-92a}, Merle presented a continuation of the explicit
blowup solution~$\psiexalpha$ beyond the singularity. To do
that, he considered the solution~$\psi^{(\varepsilon)}(t,r)$ of the
critical NLS~(\ref{eq:dDim_CriticalNLS}) with the initial condition
\begin{equation}\label{eq:MerleIC}
\psi_0^{(\varepsilon)}(r)=(1-\varepsilon)\psiexalpha(t=0,r),\qquad
0<\varepsilon\ll1.
\end{equation}
Since the power of~$\psi^{(\varepsilon)}$ is below~$\Pcr$, it exists
globally. Therefore, it is possible to continue the singular
solution~$\psiexalpha$ beyond the singularity, by taking the limit
of~$\psi^{(\varepsilon)}$ as~$\varepsilon\rightarrow0+$.
The limiting solution is given in the following Theorem:
\begin{thrm}\label{Theorem:MerleMainTheorem}
\cite{Merle-92a} Let~$\psi^{(\varepsilon)}(t,r)$ be the solution of
the critical NLS~(\ref{eq:dDim_CriticalNLS}) with the initial
condition~(\ref{eq:MerleIC}). Then, for any~$\theta\in\mathbb{R}$,
there exists a sequence $\varepsilon_n\rightarrow 0+$ (depending
on~$\theta$), such that
\begin{equation}\label{eq:MerleWeakSolution}
\psi^{(\varepsilon_n)}(t,r)\stackrel{H^1}\longrightarrow\left\{
\begin{array}{l l}
  \psiexalpha(t,r) &\quad 0\leq t<T_c,\\
  e^{i\theta}\psiexalpha^\ast(2T_c-t,r) &\quad T_c<t<\infty,
\end{array} \right.
\end{equation}
where~$\psiex(t,r)$ is given
by~(\ref{eq:ExplicitBlowupSolutionAlpha}).
\end{thrm}

Theorem~\ref{Theorem:MerleMainTheorem} shows that after the
singularity, the limiting solution is completely determined, up to
multiplication by~$e^{i\theta}$, and is symmetric with respect
to~$T_c$, i.e.,
\begin{equation}
\lim_{\varepsilon_n\rightarrow0+}\psi^{(\varepsilon_n)}(T_c+t,r)=e^{i\theta}\lim_{\varepsilon_n\rightarrow0+}\psi^{\ast(\varepsilon_n)}(T_c-t,r),\qquad
t>0.
\end{equation}

\subsection{Simulations}\label{subsec:Explicit blowup solutions}

In order to provide a numerical illustration of
Theorem~\ref{Theorem:MerleMainTheorem},
let~$\psi^{(\varepsilon)}(t,x)$ be the solution of the
one-dimensional critical NLS
\begin{subequations}\label{eq:critical_NLS_1D_problem}
\begin{equation}\label{eq:critical_NLS_1D}
i\psi_t(t,x)+\psi_{xx}+|\psi|^4\psi=0,
\end{equation}
with the initial condition
\begin{equation}\label{eq:SubcriticalInitCond}
\psi_0^{(\varepsilon)}(x)=(1-\varepsilon)\psiex(0,x)=(1-\varepsilon)\frac1{T_c^{1/2}}R\left(\frac{x}{T_c}\right)e^{i\frac1{T_c}-i\frac{x^2}{4T_c}},\qquad
T_c=0.25.
\end{equation}
\end{subequations}

\begin{figure}[ht!]
\begin{center}
\scalebox{0.8}{\includegraphics{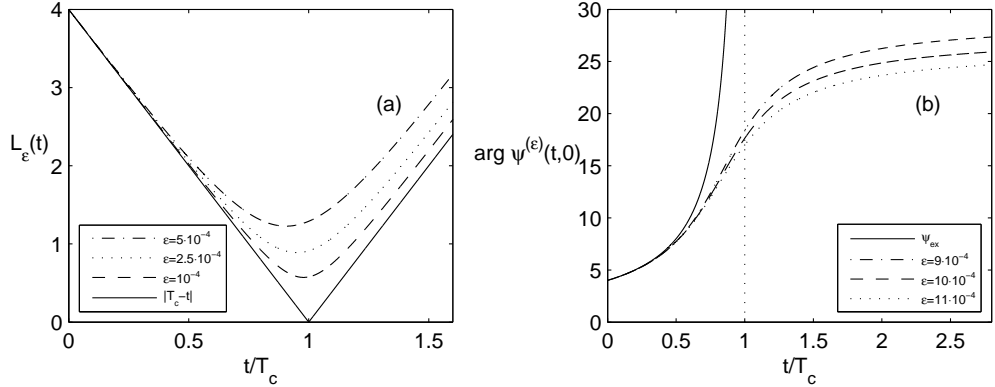}} \caption{Solution
of~(\ref{eq:critical_NLS_1D_problem}). (a):~$L_{\varepsilon}(t)$
for~$\varepsilon=5\cdot10^{-4}, 2.5\cdot10^{-4}$ and~$10^{-4}$.
Solid line is~$L=|T_c-t|$.~(b):~Accumulated phase at~$x=0$
for~$\varepsilon=0.9\cdot10^{-3}, 10^{-3}$, and~$1.1\cdot10^{-3}$.
Solid line is~${\it \arg} (\psiex(t,0))$.} \label{fig:LimExpSol_Ar}
\end{center}
\end{figure}

Let us define the width of the solution~$\psi^{(\varepsilon)}(t,x)$
as
\begin{equation}\label{eq:BeamWidthDefinition}
L_\varepsilon(t):=\left|\frac{R(0)}{\psi^{(\varepsilon)}(t,0)}\right|^{2/d},
\end{equation}
see~(\ref{eq:ExplicitBlowUpSolutions_1})
and~(\ref{eq:AdiabaticCollapse}). Figure~\ref{fig:LimExpSol_Ar}(a)
shows
that~$\lim_{\varepsilon\rightarrow0+}L_\varepsilon(t)=|T_c-t|$, both
for~$t<T_c$ and for~$T_c<t$, in agreement with
Theorem~\ref{Theorem:MerleMainTheorem}.

In order to observe the loss of the phase after the singularity, we
compute the effect of small changes in the initial condition on the
phase, by solving equation~(\ref{eq:critical_NLS_1D_problem})
with~$\varepsilon=0.9\cdot10^{-3}, 10^{-3}$ and~$1.1\cdot10^{-3}$.
Figure~\ref{fig:LimExpSol_Ar}(b) shows that
these~$\mathcal{O}(10^{-4})$ changes in the initial condition lead
to~$\mathcal{O}(1)$ changes in the phase for~$t>T_c$, which is a
manifestation of the post-collapse phase loss
as~$\varepsilon\rightarrow0+$.

\section{Sub-threshold power continuation}\label{sec:NonExplicitBeam}

The main weakness of Theorem~\ref{Theorem:MerleMainTheorem} is that
it only applies to the explicit
solutions~$\psiexalpha$.
 We now generalize Merle's
continuation to{\rev{generic}} initial profile~$F(\X)\in H^1$, as
follows. Consider the solution~$\psi(t,\X;K)$ of the critical
NLS~(\ref{eq:dDim_CriticalNLS}) with the initial condition
\begin{equation}\label{eq:GaussianIC}
\psi_0(\X;K)=K\cdot F(\X),\qquad F(\X)\in H^1,\qquad K>0.
\end{equation}
Let
\begin{equation}\label{eq:Kth_Def}
K_{th}=\inf\{K~|~\psi(t,\X;K)\ \mbox{collapses at
some~$0<T_c(K)<\infty$}\}.
\end{equation}
Let us consider the case where the infimum is attained, i.e.,
when~$\psi(t;\X;K_{th})$ becomes singular at a finite time. By
construction, the solution~$\psi^{(\varepsilon)}$ of the critical
NLS~(\ref{eq:dDim_CriticalNLS}) with the initial condition
\begin{equation}\label{eq:GeneralInitialCondition}
\psi_0^{(\varepsilon)}(\X)=(1-\varepsilon)\cdot
\psi_0^{(F)}(\X),\qquad\psi_0^{(F)}(\X)=K_{th}\cdot F(\X),
\end{equation}
exists globally for~$0<\varepsilon\ll1$, but collapses
for~$-1\ll\varepsilon<0$. Therefore, as in
Theorem~\ref{Theorem:MerleMainTheorem}, we can define the
continuation of~$\psi(t,\X;K_{th})$ beyond the singularity, by
considering the limit of~$\psi^{(\varepsilon)}$
as~$\varepsilon\rightarrow0+$. Using asymptotic analysis, in
Section~\ref{subsec:ProofOfProp1} we derive the following result,
which is the non-rigorous asymptotic analog of
Theorem~\ref{Theorem:MerleMainTheorem}:
\begin{prop}\label{Proposition:GaussianLimitSolution}
Let~$\psi^{(\varepsilon)}(t,r)$ be the solution of the critical
NLS~(\ref{eq:dDim_CriticalNLS}) with the initial
condition~(\ref{eq:GeneralInitialCondition}), where~$F(r)$ is
radial. Assume that~$\psi^{(\varepsilon=0)}(t,r)$ becomes singular
at~$0<T_c<\infty$. Then, for any~$\theta\in\mathbb{R}$, there exists
a sequence~$\varepsilon_n\rightarrow0+$ (depending on~$\theta$)
and a function $\phi \in L^2$,
such that
\begin{equation}\label{eq:AsymptoticProfileGaussian}
\begin{array}{l l}
\lim_{t\rightarrow
T_c-}\left[\lim_{\varepsilon\rightarrow0+}\psi^{(\varepsilon)}(t,r)-\psiexalpha(t,r)e^{i\theta_0}\right]=\phi(r)
&\\
\quad=\lim_{t\rightarrow
T_c+}\left[\lim_{\varepsilon_n\rightarrow0+}\psi^{(\varepsilon_n)}(t,r)-\psiexalpha^*(2T_c-t,r)e^{i\theta}\right],
\end{array}
\end{equation}
where the above limits are in~$L^2$,~$\psiexalpha$ is given
by~(\ref{eq:ExplicitBlowupSolutionAlpha}),~$T_c$,~$\alpha\in\mathbb{R}^+$
and~$\theta_0\in\mathbb{R}$. Therefore, locally near the
singularity, the limiting solution satisfies
Properties~1 and~2.

In particular, the limiting width of the collapsing core is given by
\begin{equation}
   \label{item:GaussianPhaseNonUniqueness}
L^{(0)}(t):=\lim_{\varepsilon\rightarrow0+}L(t;\varepsilon)\sim\alpha|t-T_c|,\qquad
t\rightarrow T_c\pm.
\end{equation}
\end{prop}

By Theorem~\ref{Theorem:minimal_power}, if $\psi(t;\X;K_{th})$ is a singular solution which is not given by~\eqref{eq:minimal_power_generic_general},
then $||\psi(t;\X;K_{th})||_2^2 >\Pcr$.
In that case,
$$
||\phi||_2^2 = \lim_{\varepsilon\rightarrow0+}||\psi^{(\varepsilon)}||_2^2-\Pcr
= ||\psi(t;\X;K_{th})||_2^2-\Pcr>0 .
$$

\begin{corol}
If
$$
\psi_0^{(F)}(\X) \not= \frac{1}{\alpha T_c} R^{(0)} \left(
\frac{|{\bf x}-{\bf x}_0|}{\alpha T_c} \right) e^{i\theta
  -i \frac{|{\bf x}-{\bf x}_0|^2}{4 T_c}
+ i \frac{{\bf c} \cdot ({\bf x}-{\bf x}_0)}{2}},
$$
 the limiting solution in Continuation Result~\ref{Proposition:GaussianLimitSolution} is a Bourgain-Wang solution~$\psi_{\rm BW}$, both before and after the singularity.
\end{corol}

\subsection{Simulations}
\label{subsec:CriticalGaussianBeamNumerical}

\begin{figure}[ht!]
\begin{center}
\scalebox{0.7}{\includegraphics{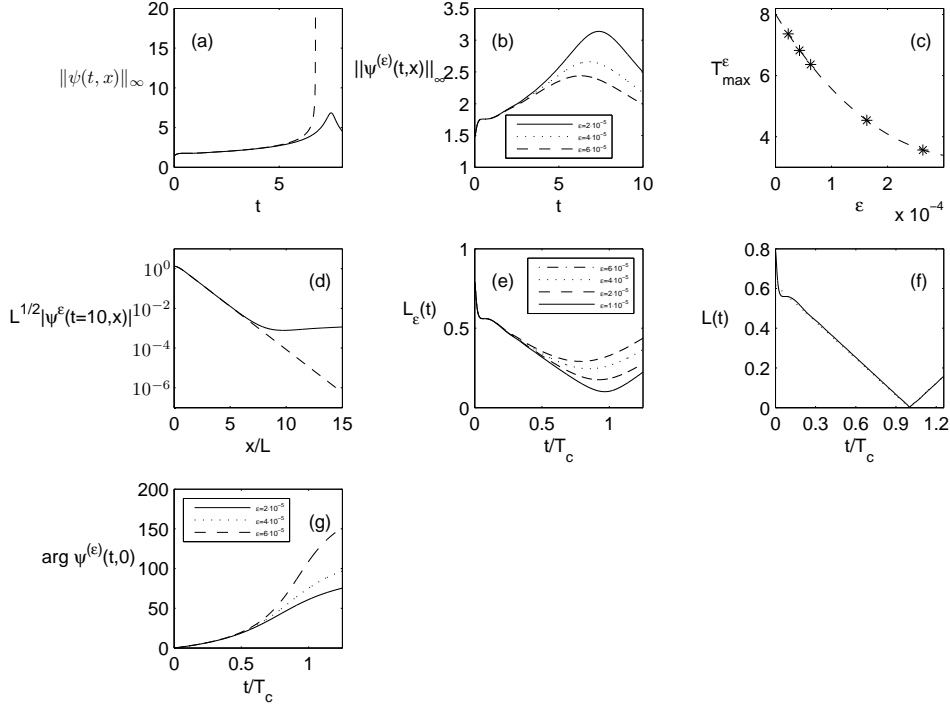}}
\caption{Solution of the one-dimensional critical
NLS~(\ref{eq:critical_NLS_1D}) with: (a)~The initial
condition~(\ref{eq:GaussianIC}) with~$K=1.48139$ (solid)
and~$K=1.48140$ (dashes).~(b)--(g)~the initial
condition~(\ref{eq:GenericSubcriticalInitCond})
with~$K_{th}=1.481395$ with various values
of~$\varepsilon$.~(b)~Solution for~$\varepsilon=2\cdot10^{-5},
4\cdot10^{-5}$, and~$6\cdot10^{-5}$. (c)~The maximum focusing
time~$T_{max}^{\varepsilon}$. (d)~Rescaled solution
for~$\varepsilon=2\cdot10^{-5}$ at~$t=10\approx1.35
T_{max}^\varepsilon$ (solid), and the~$R(x)$ profile (dashes), on a
semi-logarithmic scale. (e)~$L_\varepsilon(t)$
for~$\varepsilon=6\cdot10^{-5},4\cdot10^{-5},2\cdot10^{-5}$
and~$1\cdot10^{-5}$. (f)~Extrapolation of~$\{L_\varepsilon(t)\}$
from~(e) to~$\varepsilon=0$ (solid). Dashed line is
$L=\alpha|T_c-t|$ with $\alpha =  0.0777$ and $T_c = 7.96$.
(g)~Accumulated phase at~$x=0$.}\label{fig:MelreGeneralization}
\end{center}
\end{figure}

In order to illustrate the results of
Continuation Result~\ref{Proposition:GaussianLimitSolution}, we solve
numerically the one-dimensional critical
NLS~(\ref{eq:critical_NLS_1D}) with the initial condition
\begin{equation}\label{eq:GenericSubcriticalInitCond}
\psi_0^{(\varepsilon)}(x)=(1-\varepsilon)\cdot\psi_0^{(F)}(x),\qquad\psi_0^{(F)}(x)=K_{th}e^{-x^2}.
\end{equation}
We first compute the value of~$K_{th}$. We solve the
NLS~(\ref{eq:critical_NLS_1D}) with the initial
condition~$\psi_0=Ke^{-x^2}$.
Figure~\ref{fig:MelreGeneralization}(a) shows that for~$K=1.48140$
the beam collapses, while for~$K=1.48139$ the collapse is arrested.
Therefore,~$K_{th}=1.481395\pm5\cdot10^{-6}$.

In Figure~\ref{fig:MelreGeneralization}(b) we plot the solution
of~(\ref{eq:critical_NLS_1D}) with the initial
condition~(\ref{eq:GenericSubcriticalInitCond}), as a function
of~$t$. The solution amplitude increases up
to~$t=T_{max}^{\varepsilon}:=\argmax_t\|\psi^{(\varepsilon)}(t,x=0)\|_\infty$,
and then decreases. As expected, both~$T_{max}^{\varepsilon}$ and
the maximal
amplitude~$|\psi^{(\varepsilon)}(T^{\varepsilon}_{max},0)|$ increase
as~$\varepsilon\rightarrow0+$. In
Figure~\ref{fig:MelreGeneralization}(c) we
plot~$T_{max}^{\varepsilon}$ as a function of~$\varepsilon$.
Extrapolation of these values to~$\varepsilon=0$ shows that
\begin{equation}\label{eq:Tc_Range}
T_c:=\lim_{\varepsilon\rightarrow0+}T_{max}^\varepsilon\approx8.00.
\end{equation}
Since~$T_c<\infty$,~$\psi(t,x;K_{th})$ is a singular solution.

In order to confirm that the profile of the collapsing core is a
rescaled~$R^{(0)}$ profile, see
equation~(\ref{eq:AsymptoticProfileGaussian}), in
Figure~\ref{fig:MelreGeneralization}(d) we
plot~$|\psi^{(\varepsilon)}|$ for~$\varepsilon=2\cdot10^{-5}$ at
~$t=10\approx1.35T_{max}^\varepsilon$, i.e., after it collapse has
been arrested, and observe that for~$0\leq x/L(t)\leq5$, the
rescaled profile is indistinguishable from~$R^{(0)}$, while
for~$6<x/L(t)$ the two curves are different. This confirms that the
inner core collapses with the~$\psiexalpha$ profile, but the outer
tail does not.

In Figure~\ref{fig:MelreGeneralization}(e) we plot the solution
width~$L_\varepsilon(t)$ for various values of~$\varepsilon$. The
extrapolation of these curves to~$\varepsilon=0$ is in good
agreement with the predicted linear
limit~\eqref{item:GaussianPhaseNonUniqueness} for $0.25 \le  t/T_c
\le 1.2$, see Figure~\ref{fig:MelreGeneralization}(f). Finally, in
order to observe the loss of the phase after the singularity, we
compute the effect of~$\mathcal{O}(10^{-5})$ changes in the initial
condition on the phase, by solving the one-dimensional
NLS~(\ref{eq:critical_NLS_1D}) with the initial
condition~(\ref{eq:GenericSubcriticalInitCond})
with~$\varepsilon=2\cdot10^{-5}, 4\cdot10^{-5}$ and~$6\cdot10^{-5}$.
Figure~\ref{fig:MelreGeneralization}(g) shows that
these~$\mathcal{O}(10^{-5})$ changes in the initial condition lead
to~$\mathcal{O}(1)$ changes in the phase for~$8\leq t\leq10$. Since
~$T_c\approx8$, see~(\ref{eq:Tc_Range}), these~$\mathcal{O}(1)$
changes occur after the singularity, in agreement with
Continuation Result~\ref{Proposition:GaussianLimitSolution}.

\subsection{Proof of
Continuation Result~\ref{Proposition:GaussianLimitSolution}}
  \label{subsec:ProofOfProp1}

\subsubsection{Adiabatic collapse -
review}\label{subsec:AsymptoticAnalysis_GeneralIC}

In order to
compute asymptotically the limit of~$\psi^{\varepsilon}$
as~$\varepsilon\rightarrow0+$, we recall that, in general, the
collapse of radial solutions with power close to~$\Pcr$ can be
divided into two stages, see~\cite{PNLS-99}:
\begin{enumerate}
\item During the initial non-adiabatic self-focusing stage, the
solution "splits" into a collapsing core~$\psi_{\rm core}$ and a
non-collapsing "tail", i.e.,
\begin{subequations}
   \label{eq:AdiabaticCollapse}
\begin{equation}
   \label{eq:AdiabaticCollapse-1}
\psi(t,r)\sim\left\{
\begin{array}{l l}
  \psi_{\rm core} & \quad 0\leq\frac{r}{L(t)}\leq\rho_0,\\
  \psi_{\rm tail} & \quad\quad \rho_0\ll\frac{r}{L(t)},
\end{array} \right.
\end{equation}
where~$L(t)$ is the collapsing core "width"
and~$\rho_0=\mathcal{O}(1)$.
\item As~$t\rightarrow T_c$,~$\psi_{\rm core}$ approaches the self-similar profile
\begin{equation}
   \label{eq:AsymptoticProfile}
\psi_{R^{(0)}}=\frac1{L^{d/2}(t)}R^{(0)}(\rho)e^{iS},\quad
\rho=\frac{r}{L},\quad S=\tau(t)+\frac{L_t}{L}\frac{r^2}{4},\quad
\tau=\int_0^t\frac1{L^2(s)}ds,
\end{equation}
\end{subequations}
where~$R^{(0)}$ is the ground state
of~(\ref{eq:dDimCritical_R_ODE}). In contrast, the "tail" continues
to propagate forward. In particular,~$\lim_{t\rightarrow
T_c}\psi_{\rm tail}=\phi(\X)\in L^2$, see
Theorem~\ref{thm:Merle+Raphael}.

Once the profile of~$\psi_{\rm core}$ is close enough
to~$\psi_{R^{(0)}}$, the dynamics of the collapsing core becomes
nearly adiabatic, and is governed, to leading order, by the reduced
equations~\cite{Papanicolaou-88_1,Papanicolaou-88_2,Fraiman-1985}
\begin{equation}\label{eq:ReducedEqGaussian}
\beta_t(t)=-\frac{\upsilon(\beta)}{L^2},\qquad
L_{tt}=-\frac{\beta}{L^3},
\end{equation}
where
\begin{subequations}\label{eq:ReducedEqCont}
\begin{equation}\label{eq:ni_beta_def}
\upsilon(\beta)=\left\{
\begin{array}{l l}
  c_{\nu}e^{-\pi/\sqrt{\beta}}, & \qquad \beta>0,\\
  0, & \qquad \beta\leq0,\\
\end{array} \right.
\end{equation}
and
\begin{equation}\label{eq:ModulationTheoryDefs}
c_{\nu}=\frac{2A^2_R}{M},\quad
A_R=\lim_{r\rightarrow\infty}e^rr^{(d-1)/2}R^{(0)}(r),\quad
M=\frac1{4}\int_0^\infty r^2(R^{(0)})^2r^{d-1}dr.
\end{equation}
\end{subequations} In addition, the parameter~$\beta$ is proportional to the
excess power above~$\Pcr$ of the collapsing core~$\psi_{\rm core}$,
i.e.,
\begin{equation}\label{eq:BetaDefGaussian}
\beta\sim\frac{\|\psi_{\rm core} \|_2^2-\Pcr}{M},
\end{equation}
see~\cite{Malkin-1993,PNLS-99}.
\end{enumerate}

\subsubsection{Asymptotic analysis} For simplicity, we assume that
the initial condition~(\ref{eq:GeneralInitialCondition}) is radial.
By construction, for the initial
condition~(\ref{eq:GeneralInitialCondition}),
\[
\left\{
\begin{array}{l l}
  \|\psi_{\rm core}\|_2^2>\Pcr & \quad \varepsilon<0,\\
  \|\psi_{\rm core}\|_2^2<\Pcr & \quad \varepsilon>0.\\
\end{array} \right.
\]
Therefore, by~(\ref{eq:BetaDefGaussian}),
\begin{equation}\label{eq:beta_vs_epsilon}
\left\{
\begin{array}{l l}
  \beta>0 & \quad \varepsilon<0,\\
  \beta<0 & \quad \varepsilon>0.\\
\end{array} \right.
\end{equation}
Hence, by~(\ref{eq:ni_beta_def}),~$\upsilon(\beta)\equiv0$
for~$\varepsilon>0$. Therefore~$\beta_t(t;\varepsilon)\equiv0$.
Thus, the self-focusing dynamics is governed by
\begin{equation}\label{eq:ReducedEqGaussianSol}
L_{tt}(t;\varepsilon)=-\frac{\beta(\varepsilon)}{L^3},
\end{equation}
where~$\beta(\varepsilon)$ is independent of~$t$.

Let~$t_{ad}(\varepsilon)>0$ denote the time at
which~$\psi^{(\varepsilon)}$ "enters" the adiabatic stage, i.e.,
when~$\psi^{(\varepsilon)}\sim\psi_{R^{(0)}}$, so that
equations~(\ref{eq:ReducedEqGaussian}) hold. Therefore, for~$t\geq
t_{ad}(\varepsilon)$, the dynamics is given by
\begin{equation}\label{eq:L_ODE_Guassian}
L_{tt}(t)=-\frac{\beta(\varepsilon)}{L^3},\qquad
L(t_{ad})=L_{ad}(\varepsilon),\quad
L_t(t_{ad})=L'_{ad}(\varepsilon).
\end{equation}
\begin{lem}\label{Lemma:LimitReducedEqGaussianSol}
Let~$L(t;\varepsilon)$ be the solution of~(\ref{eq:L_ODE_Guassian}).
Denote
\begin{equation}\label{eq:t_ad_definition}
t_{ad}:=\lim_{\varepsilon\rightarrow0+}t_{ad}(\varepsilon),\qquad\alpha:=\lim_{\varepsilon\rightarrow0+}L_t(t_{ad}(\varepsilon);\varepsilon),\qquad
T_c:=\lim_{\varepsilon\rightarrow0+}\left(-\frac{L_{ad}(\varepsilon)L'_{ad}(\varepsilon)}{(L'_{ad}(\varepsilon))^2-\frac{\beta(\varepsilon)}{L_{ad}^2(\varepsilon)}}\right).
\end{equation}
Then,
\begin{equation}\label{eq:GaussianLimitBeamWidth}
\lim_{\varepsilon\rightarrow0+}L(t;\varepsilon)=\alpha|t-T_c|,\qquad
t_{ad}\leq t<\infty.
\end{equation}
\end{lem}
\Beginproof Two integrations show that the solution
of~(\ref{eq:L_ODE_Guassian}) is given by
\[
L^2(t;\varepsilon)=c_1(\varepsilon)\left(t+\frac{L_{ad}(\varepsilon)L'_{ad}(\varepsilon)}{c_1(\varepsilon)}\right)^2-\frac{\beta(\varepsilon)}{c_1(\varepsilon)},
\qquad
c_1(\varepsilon)=(L'_{ad}(\varepsilon))^2-\frac{\beta(\varepsilon)}{L_{ad}^2(\varepsilon)}.
\]
By~(\ref{eq:beta_vs_epsilon}),~$\lim_{\varepsilon\rightarrow0+}\beta(\varepsilon)=0$.
Therefore,~$\lim_{\varepsilon\rightarrow0+}c_1(\varepsilon)=\alpha^2$.
Hence,~$\lim_{\varepsilon\rightarrow0+}L^2(t;\varepsilon)=\alpha^2(t-T_c)^2$.~\Endproof
\begin{corol}\label{Corollary:LimitPhaseGaussianTc}
$\lim_{\varepsilon\rightarrow0+}\tau(t=T_c;\varepsilon)=\infty$.
\end{corol}
\Beginproof By~(\ref{eq:AsymptoticProfile})
and~(\ref{eq:GaussianLimitBeamWidth}),
\[
\lim_{\varepsilon\rightarrow0+}\tau(t=T_c;\varepsilon)=\lim_{\varepsilon\rightarrow0+}\int_0^{T_c}\frac{dt}{L^2(t;\varepsilon)}=
\int_0^{T_c}\frac{dt}{\alpha^2(T_c-t)^2}=\infty.
\]\Endproof

In order to go back from~$L(t)$ to~$\psi$, let us note that:
\begin{lem}\label{Lemma:AsymptoticProfile_L}
Let $\psi_{R^{(0)}}$ be given by~\eqref{eq:AsymptoticProfile}. If~$L(t)=\alpha|T_c-t|$, then
\begin{equation}
\psi_{R^{(0)}}(t,r)=\left\{
\begin{array}{l l}
  e^{i \theta_0} \psiexalpha(t,r),&\qquad
0\leq t<T_c,\\
e^{i \theta_1} \psiexalpha^*(2T_c-t,r),&\qquad t>T_c,
\end{array} \right.
\end{equation}
where~$\psiexalpha(t,r)$ is given
by~(\ref{eq:ExplicitBlowupSolutionAlpha}), and $\theta_0, \theta_1 \in {\Bbb R}$.
\end{lem}

\subsubsection{Proof of
Continuation Result~\ref{Proposition:GaussianLimitSolution}}

In Lemma~\ref{Lemma:LimitReducedEqGaussianSol} we saw that the
solution~$L(t)$ of the reduced system is given
by~(\ref{eq:GaussianLimitBeamWidth}) for~$t_{ad} \le t<\infty$.
Therefore, by Lemma~\ref{Lemma:AsymptoticProfile_L},
when~$t_{ad}\leq t<T_c$,~$\psi_{R^{(0)}}(t,r)=  \psiexalpha(t,r) e^{
i \theta_0}$, and when~$T_c<t$,~$\psi_{R^{(0)}}(t,r)=
\psiexalpha^\ast(2T_c-t,r) \lim_{\epsilon_n \to 0}e^{ i
\theta_1(\epsilon_n)}$.

Since~$\arg\psi(t,0)\sim\arg\psi_{R^{(0)}}(t,0)=\tau(t)$,
Corollary~\ref{Corollary:LimitPhaseGaussianTc} shows that the
limiting phase becomes infinite at~$T_c$, hence also for~$t>T_c$.
Therefore, for a given~$t>T_c$ and~$\theta\in\mathbb{R}$, there
exists a sequence~$\varepsilon_n\rightarrow0+$, such
that~$\lim_{\varepsilon_n\rightarrow0+}\arg\psi^{(\varepsilon_n)}(t,0)=\theta$.
Since as~$t\rightarrow T_c$,~$\psi_{\rm core}\rightarrow\psi_{R^{(0)}}$
and~$\psi_{\rm tail}\rightarrow\phi(\X)\in L^2$, the Proposition
follows.

\section{Time-reversible continuations}
 \label{sec:SymmetryProperty}

The continuation in Continuation
Result~\ref{Proposition:GaussianLimitSolution} preserves
Properties~1 and~2 of
Merle's first continuation. We now show that these two properties
hold for continuations of the NLS that preserve the NLS invariance
under the transformation,
\begin{equation}
   \label{eq:NLSInvariantTransformations}
t\rightarrow-t\quad\mbox{and}\quad\psi\rightarrow\psi^\ast,
\end{equation}
and also satisfy some additional conditions.
\begin{prop}
\label{Proposition:WeakSolutionSymmetry} Let~$\psi(t,\X)$ be a
solution of the NLS~(\ref{eq:DdimensionalNLS}) that blows up
at~$T_c$, and let~$\psi^{(\varepsilon)}(t,\X)$ be a smooth
continuation of~$\psi(t,\X)$, such that
\begin{enumerate}
\item $\psi^{(\varepsilon)}$ exists globally
for~$0<\varepsilon\ll 1$.
\item~$\lim_{\varepsilon\rightarrow0+}\psi^{(\varepsilon)}(t,\X)=\psi(t,\X)$ in~$L^{2 \sigma + 2}$
for~$0\leq t<T_c$.
\item~$\psi^{(\varepsilon)}$ is invariant under the
transformation~(\ref{eq:NLSInvariantTransformations}).
\item $\lim_{t\rightarrow T_c}\arg\psi(t,0)=\infty$.
\item Let~$T^{\varepsilon}_{max}:=\argmax_t\|\psi^{(\varepsilon)}(t,r)\|_{2 \sigma + 2}$
denote the time at which the collapse of~$\psi^{(\varepsilon)}$ is
arrested. Then, for all ${\bf x} \in {\Bbb R}^d$,
\begin{equation}
\label{eq:cond-collimated}
\arg\psi^{(\varepsilon)}(T^{\varepsilon}_{max},\X)\equiv\alpha^{(\varepsilon)},\qquad\alpha^\varepsilon\in\mathbb{R}.
\end{equation}
\end{enumerate}
Then, for any~$\theta\in\mathbb{R}$, there exists a
sequence~$\varepsilon_n\rightarrow0+$ (depending on~$\theta$), such
that
\[
\lim_{\varepsilon_n\rightarrow0+}\psi^{(\varepsilon_n)}(T_c+t,\X)=e^{i\theta}\psi^\ast(T_c-t,\X),\qquad
t>0.
\]
Hence, the continuation satisfies Properties~1
and~2.
\end{prop}
\Beginproof Assume first that~$\alpha^{(\varepsilon)}=0$. Then, from
the invariance of~$\psi^{(\varepsilon)}$
under~(\ref{eq:NLSInvariantTransformations}) it follows that
\begin{equation}\label{eq:PsiEpsilonSymmetry_1}
\psi^{(\varepsilon)}(T_{max}^\varepsilon+t,\X)=\psi^{\ast(\varepsilon)}(T_{max}^\varepsilon-t,\X),\qquad
t>0.
\end{equation}
If~$\alpha^{(\varepsilon)}\neq0$,
then~(\ref{eq:PsiEpsilonSymmetry_1}) holds
for~$e^{-i\alpha^{(\varepsilon)}}\psi^{(\varepsilon)}$, i.e.,
\begin{equation}\label{eq:NLS_SymmetryRes}
e^{-i\alpha^{(\varepsilon)}}\psi^{(\varepsilon)}(T^{\varepsilon}_{max}+t,\X)=e^{i\alpha^{(\varepsilon)}}\psi^{\ast(\varepsilon)}(T^{\varepsilon}_{max}-t,\X),\qquad
t>0.
\end{equation}
Hence,
\begin{equation}\label{eq:PsiEpsilonSymmetry_2}
\psi^{(\varepsilon)}(T^{\varepsilon}_{max}+t,\X)=e^{i\theta^{(\varepsilon)}}\psi^{\ast(\varepsilon)}(T^{\varepsilon}_{max}-t,\X),\qquad
t>0,
\end{equation}
where~$\theta^{(\varepsilon)}=2\alpha^{(\varepsilon)}$.
Since~$\lim_{\varepsilon\rightarrow0+}\theta^{(\varepsilon)}=\lim_{\varepsilon\rightarrow0+}\arg\psi^{(\varepsilon)}(T_{max}^\varepsilon,0)=\infty$,
there exists~$\varepsilon_n\rightarrow0+$ such
that~$\theta=\lim_{\varepsilon_n\rightarrow0+}\theta^{(\varepsilon_n)}$.
In addition, from~$\lim_{t \to T_c} ||\psi (t)||_{2 \sigma +2} =
\infty$ and~$\lim_{\epsilon \to 0} ||\psi^\epsilon(t)||_{2 \sigma
+2} = ||\psi (t)||_{2 \sigma +2}$ for $0 \le t<T_c$, it follows
that~$\lim_{\varepsilon\rightarrow0+}T_{max}^{\varepsilon}=T_c$. The
result follows by taking the limit
of~(\ref{eq:PsiEpsilonSymmetry_2}). \Endproof

Continuation Result~\ref{Proposition:WeakSolutionSymmetry} holds for both
the critical and the supercritical NLS. The interpretation of the
conditions of Continuation Result~\ref{Proposition:WeakSolutionSymmetry} is
as follows. Conditions~1 and~2 say that $\psi^{(\epsilon)}$ is a
continuation of~$\psi$. Condition~3 says that the continuation is
time reversible. Condition~4 says that the phase of the singular
solution becomes infinite at the singularity. This condition holds
for all known singular solutions of the critical and supercritical
NLS. Condition~5 says that at~$t=T^{\varepsilon}_{max}$, the
solution is collimated. Intuitively, this  is because the solution
is focusing for $t<T^{\varepsilon}_{max}$, and defocusing for
$t>T^{\varepsilon}_{max}$.

We now confirm that the sub-threshold power continuation of
Continuation Result~\ref{Proposition:GaussianLimitSolution} satisfies the conditions
of Continuation Result~\ref{Proposition:WeakSolutionSymmetry}.
By~(\ref{eq:AsymptoticProfile}),
\[
\psi^{(\varepsilon)}(t,r)\sim\frac1{L_\varepsilon^{d/2}(t)}R^{(0)}\left(\frac{r}{L}\right)e^{i\tau_\varepsilon+\frac{(L_\varepsilon)_t}{L_\varepsilon}\frac{r^2}{4}},\qquad\frac{d\tau_\varepsilon}{dt}=\frac1{L_\varepsilon^2}.
\]
Since~$L_\varepsilon(t)$ attains its minimum
at~$T_{max}^\varepsilon$,
then~$L_\varepsilon'(T_{max}^\varepsilon)=0$.
Therefore,~$\arg\psi^{(\varepsilon)}(T^\varepsilon_{max},r)\equiv\alpha^{(\varepsilon)}$,
where~$\alpha^{(\varepsilon)}=\tau_\varepsilon(T^\varepsilon_{max})$.
Furthermore,
since~$L_\varepsilon^2(t)\sim\frac1{\|\nabla\psi^{(\varepsilon)}\|_2^2} \sim \frac1{\|\nabla\psi\|_2^2}$,
and since in the critical NLS~$\|\nabla\psi\|_2^2\geq
M(T_c-t)^{-1}$~\cite{Cazenave-90},
$$
\lim_{\varepsilon\rightarrow0+} \tau_\varepsilon
(T^\varepsilon_{max}) =  \lim_{\varepsilon\rightarrow0+}
\int_0^{T^\varepsilon_{max}} \frac{dt}{L_\varepsilon^2}  \ge M
\int_0^{T_c} \frac{dt}{T_c-t} =  \infty.
$$
Therefore,~$\lim_{\varepsilon\rightarrow0+}\alpha^{(\varepsilon)}=\infty$.
Hence,~$\lim_{\varepsilon\rightarrow0+}\theta^{(\varepsilon)}=\infty$.

From all the conditions of
Continuation Result~\ref{Proposition:WeakSolutionSymmetry}, the only one
whose validity is questionable is condition~5. It is reasonable to
expect that this condition would hold for the collapsing core. There
is no reason, however, why it should hold for the non-collapsing
tail.

An immediate consequence of Continuation Result~\ref{Proposition:WeakSolutionSymmetry} is that
\begin{corol}
   \label{corol:point-sing}
Under the conditions of
Continuation Result~\ref{Proposition:WeakSolutionSymmetry}, the limiting
solution $\lim_{\varepsilon_n \rightarrow 0+} \psi^{(\varepsilon_n)}$ is
in~$H^1 $for $t>T_c$. Hence,  the continuation leads to a point
singularity, and not to a filament singularity.
\end{corol}

In Section~\ref{sec:nl-saturation} we will see that time-reversible
continuations can also lead to a filament singularity. In that case,
however, condition~5 does not hold.

\section{Vanishing nonlinear-saturation continuation}
\label{sec:nl-saturation}

\subsection{Merle's second continuation}

In \cite{Merle-92b}, Merle presented a different continuation, which
is based on arresting the collapse with an addition of nonlinear
saturation.

\begin{thrm}
\label{thm:Merle-nl-saturation}{\rev{\cite{Merle-92b}}} Let
$d\geq2$, and consider radial initial data~$\psi_0(r)\in
H^1\bigcap\{r\psi_0\in L^2\}$, such that the solution~$\psi(t,\X)$
of the critical NLS~(\ref{eq:dDim_CriticalNLS}) blows up in finite
time~$T_c$. For~$\varepsilon>0$ and~$1+4/d<q<(d+2)/(d-2)$,
let~$\psi_\varepsilon(t,\X)$ be the solution of the saturated
critical NLS
\begin{equation}
 \label{SaturatedNLS}
i\psi_t(t,\X)+\Delta\psi+|\psi|^{4/d}\psi-\varepsilon|\psi|^{q-1}\psi=0,\qquad
\psi(0,\X)=\psi_0(r).
\end{equation}

If for~$T_0>T_c$, there is a constant~$C>0$ such
that~$\int|\X|^2|\psi_\varepsilon(T_0,\X)|^2 d\X\leq C$, then:
\begin{enumerate}
\item There is
a function~$\tilde{\psi}(t,\X)$ defined for $t<T_0$, such that for
all~$r_0>0$,~$\tilde{\psi}\in\mathcal{C}([0,T_0),L^2(|\X|\geq
r_0))$, and~$\psi_\varepsilon(t,\X)\rightarrow\tilde{\psi}(t,\X)$
in~$\mathcal{C}([0,T_0),L^2(|\X|\geq r_0))$
as~$\varepsilon\rightarrow0$.
\item For~$t<T_0$, there is~$m(t)\geq0$
such that~$|\psi_\varepsilon(t,\X)|^2\rightarrow
m(t)\delta(\X)+|\tilde{\psi}(t,\X)|^2$
as~$\varepsilon\rightarrow0$ in the distribution sense. Furthermore,
\begin{enumerate}

\item If~$m(t)\neq0$,
    then~$\|\psi_\varepsilon(t,\X)\|_{H^1}\rightarrow+\infty$
as~$\varepsilon\rightarrow0$ and~$m(t)\geq P_{cr}$.

\item If~$m(t)=0$, there is a constant~$c>0$ such that for
all~$\varepsilon$,~$\|\psi_\varepsilon(t,\X)\|_{H^1}<c$,
and~$\psi_\varepsilon(t,\X)\rightarrow\tilde{\psi}(t,\X)$ in~$L^2$.
\end{enumerate}
\item For all~$t<T_0$,~$m(t)+\int|\tilde{\psi}(t,\X)|^2d\X=\int|\psi_0(\X)|^2d\X$.
\end{enumerate}
\end{thrm}

Theorem~\ref{thm:Merle-nl-saturation} shows that the vanishing
nonlinear saturation continuation can lead to a filament
singularity. The condition $\int|\X|^2|\psi_\varepsilon(T_0,\X) \,
d\X\leq C$ is believed to hold generically.

\subsection{Malkin's analysis}
  \label{sec:Malkin}

In~\cite{Malkin-1993}, Malkin analyzed asymptotically the solutions
of the saturated critical NLS~\eqref{SaturatedNLS} with~$d=2$
and~$q=5$. Malkin showed that initially, the solution follows the
non-saturated NLS solution and self-focuses. Then, the collapse is
arrested by the nonlinear saturation, leading to focusing-defocusing
oscillations. During each oscillation, the collapsing core loses
(radiates) some power. As a result, the magnitude of the
oscillations decreases, so that ultimately, the solution approaches
a standing-wave solution of the saturated NLS.

\begin{figure}[ht!]
\begin{center}
\scalebox{0.7}{\includegraphics{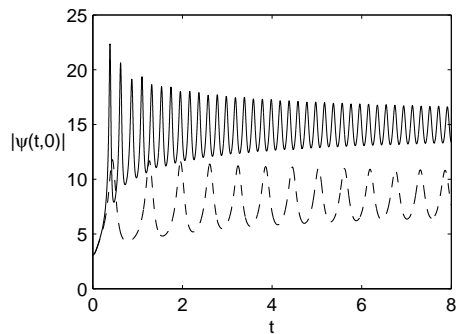}}
\caption{Solution of the saturated NLS~\eqref{SaturatedNLS}
with~$d=2$ and~$q=5$ for~$\varepsilon=0.5\cdot10^{-3}$ (solid),
and~$\varepsilon=2\cdot10^{-3}$ (dashes).
Here,~$\psi_0(r)=3.079e^{-r^2}$.}\label{fig:SaturatedNLS_Oscillations}
\end{center}
\end{figure}

If we fix the initial condition~$\psi_0$ and let $\varepsilon
\longrightarrow 0+$, then as~$\varepsilon$ decreases, the collapse
is arrested at a later stage. For example, in
Figure~\ref{fig:SaturatedNLS_Oscillations} we plot the solution of
the saturated critical NLS~\eqref{SaturatedNLS} with~$d=2$ and~$q=5$
with the initial condition~$\psi_0(r)=3.079e^{-r^2}$, and observe
that as~$\varepsilon$ decreases, the collapse is arrested at a later
stage, and the oscillations occurs at higher amplitudes. In
addition, we observe that as~$t$ increases, the oscillations
decreases. Hence, it is reasonable to assume that the amplitude of
the limiting standing-wave increases as~$\varepsilon$ decreases, and
goes to infinity as $\varepsilon \longrightarrow 0+$. In addition,
as~$\varepsilon\longrightarrow0+$, the power of the standing-wave
of~\eqref{SaturatedNLS} approaches~$\Pcr$. Therefore, Malkin's
analysis suggests that
$$
m(t) \equiv \Pcr, \qquad  T_c \le t<\infty,
$$
i.e., after the singularity the limiting solution consists of a semi-infinite filament with power~$\Pcr$,
and a regular part with power $||\psi_0||_2^2-\Pcr$.

\subsection{Importance of power radiation}

  The above results of Merle and Malkin strongly suggest that the continuation of singular NLS solutions with a vanishing nonlinear-saturation generically
leads to a filament singularity. Since the NLS with a nonlinear
saturation is time reversible, these results seem to be in
contradiction with
Continuation Result~\ref{Proposition:WeakSolutionSymmetry}, see
Corollary~\ref{corol:point-sing}. Note, however, that in
Continuation Result~\ref{Proposition:WeakSolutionSymmetry} we assumed that
the solution phase is constant at the time~$T_{max}^\epsilon$ where
its collapse is arrested (Condition~5). If this condition were to
hold for the solution of the saturated NLS, then by
Continuation Result~\ref{Proposition:WeakSolutionSymmetry}, the solution
at~$t = 2 T_{max}^{\varepsilon}$ would be given by~$\psi_0^*$. As
Figure~\ref{fig:SaturatedNLS_Oscillations} shows, this is not the
case.

The mechanism which enables the filamentation is the loss
(radiation) of power from the collapsing core to the surrounding
background. Indeed, in the absence of radiation, the collapsing core
of the solution of saturated NLS undergoes periodic oscillations,
rather than approaches a standing wave~\cite{Malkin-1993}. We stress
that the constant phase condition does hold asymptotically for the
collapsing core~\cite{Malkin-1993}. It does not, however, hold for
the regular part of the solution (the "tail").\footnote{In the
terminology of Theorem~\ref{thm:Merle+Raphael}, it holds
for~$\psi_{R^{(0)}}$, but not for~$\phi$.}

\section{Chaotic interactions}\label{sec:InteractionTwoBeams}

In Sections~\ref{sec:WeakSolNLS}--\ref{sec:SymmetryProperty} we saw
that time-reversible continuations have the property whereby the
phase becomes non-unique after the blowup time. A-priori, this phase
loss should have no effect, since multiplying the NLS solution
by~$e^{i\theta}$ does not affect the dynamics. Nevertheless, we now
show that this phase loss can affect the interaction between two
post-collapse beams (filaments).

Consider first the initial condition
\begin{equation}\label{eq:PsiTilt}
\psi_0({\bf x}) =  \psiex^{tilt,\pm{\bf{c}}}(0,\X\mp\X_0)=\psiex(0,\X\mp
{\bf{x}}_0)e^{\pm i{\bf{c}}\cdot(\X\mp
{\bf{x}}_0)/2}.
\end{equation}
By~(\ref{eq:GeneralExplicitBlowupSol}), the solution of the critical
NLS~(\ref{eq:dDim_CriticalNLS}) with the initial
condition~(\ref{eq:PsiTilt}) is given by
\[
\psiex^{tilt,\pm{\bf{c}}}(t,\X\mp\X_0)=\psiex(t,\X\mp\X_0\pm{\bf{c}}\cdot
t)e^{\pm i{\bf{c}}\cdot(\X\mp\X_0)/2-i|{\bf{c}}|^2t/4}.
\]
Therefore,~$\psiex^{tilt,\pm{\bf{c}}}(t,\X)$ is the explicit
blowup solution~(\ref{eq:ExplicitBlowupSolution}), centered
initially at~$\X=\pm\X_0$, and tilted at the angle
of~$\pm\arctan(|{\bf c}|)$.

\begin{figure}[ht!]
\begin{center}
\scalebox{0.6}{\includegraphics{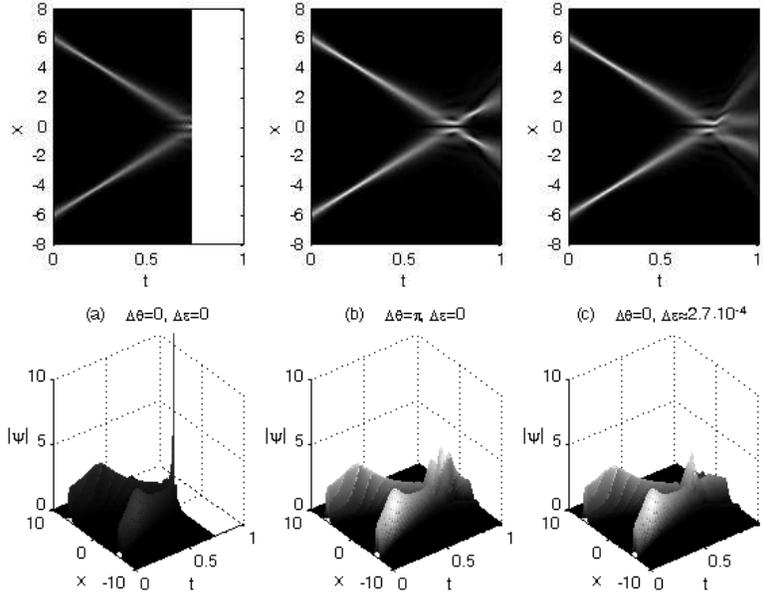}}
\caption{Solution of the one-dimensional critical
NLS~(\ref{eq:critical_NLS_1D}) with the initial
condition~(\ref{eq:DualBeamAtAnAngle})
with~$T_c=0.25$,~$x_0=6$,~$c=8$ and~$\varepsilon=10^{-3}$. (a):
In-phase identical beams~($\Delta\theta=0$
and~$\Delta\varepsilon=0$). (b): Out-of-phase identical
beams~($\Delta\theta=\pi$ and~$\Delta\varepsilon=0$). (c): In-phase
non-identical beams
~($\Delta\theta=0$~and~$\Delta\varepsilon\approx2.7\cdot10^{-4}$).}\label{fig:TwoBeamInteraction}
\end{center}
\end{figure}

We now consider the one-dimensional critical
NLS~(\ref{eq:critical_NLS_1D}) with the two tilted-beams initial
condition
\begin{equation}\label{eq:DualBeamAtAnAngle}
\psi_0(x)=(1-\varepsilon)\psiex^{tilt,+c}(0,x-x_0)+(1-(\varepsilon+\Delta\varepsilon))\psiex^{tilt,-c}(0,x+x_0))e^{i\Delta\theta},
\end{equation}
where~$\psi_{0,ex}^{tilt,\pm{\bf{c}}}$ is defined
in~(\ref{eq:PsiTilt}) and~$\varepsilon=10^{-3}$. This initial
condition correspond to two input beams, centered at $\pm x_0$,
titled toward each other, possible with a different power (when
$\Delta\varepsilon \not=0$), and with a relative phase difference
$\Delta \theta$. In Figure~\ref{fig:TwoBeamInteraction}(a) we plot
the solution when~$\Delta\varepsilon=0$ and~$\Delta\theta=0$
(equal-power, in-phase input beams). Since the power of each beam is
slightly below~$\Pcr$, each beam focuses up to a certain time
($t\approx0.3$), and then defocuses. Subsequently, the two beams
intersect around~$t\approx0.6$. Since the beams are in phase, they
interact constructively. As a result, their total power
is~$\approx2\Pcr$. Hence, the solution collapses at~$T_c\approx0.6$.

In Figure~\ref{fig:TwoBeamInteraction}(b) we repeat this simulation
with~$\Delta\varepsilon=0$ and~$\Delta\theta=\pi$ (equal-power,
out-of-phase input beams). Before the two beams intersect, their
dynamics is the same as in Figure~\ref{fig:TwoBeamInteraction}(a).
When they intersect at~$t\approx0.6$, however, the two beams are
out-of-phase. Hence, they repel each other. Since each beam has
power below~$\Pcr$, there is no collapse.

In Figure~\ref{fig:TwoBeamInteraction}(c) we repeat this simulation
with~$\Delta\theta=0$ and~$\Delta\varepsilon\approx2.7\cdot10^{-4}$
(in-phase and slightly-different input powers), and observe that
at~$t\approx0.6$ the two beams repel each other and there is no
collapse. In particular, comparison of
Figures~\ref{fig:TwoBeamInteraction}(a)
and~\ref{fig:TwoBeamInteraction}(c) shows that
\textit{the~$\mathcal{O}(10^{-4})$ change in the initial condition
lead to a completely different "post collapse" interaction pattern
between the two beams.}

\begin{figure}[ht!]
\begin{center}
\scalebox{0.5}{\includegraphics{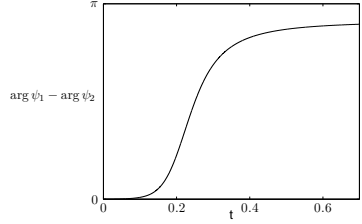}} \caption{
$\arg\psi_1(t,x=x_0-ct)-\arg\psi_2(t,x=-x_0+ct)$.}\label{fig:DeltaArgPsi_DeltaEps}
\end{center}
\end{figure}

The dynamics in Figure~\ref{fig:TwoBeamInteraction}(c) is
qualitatively the same as in Figure~\ref{fig:TwoBeamInteraction}(b).
This suggests that when the two beams in
Figure~\ref{fig:TwoBeamInteraction}(c) intersect, their phase
difference is~$\approx\pi$. Indeed, let~$\psi_1(t,x)$
and~$\psi_2(t,x)$ be the solutions of the one-dimensional critical
NLS~(\ref{eq:critical_NLS_1D}) with the initial
conditions~$(1-\varepsilon)\psiex^{tilt,+c}(0,x-x_0)$
and~$(1-(\varepsilon+\Delta\varepsilon))\psiex^{tilt,-c}(0,x+x_0))$,
respectively. In Figure~\ref{fig:DeltaArgPsi_DeltaEps} we plot the
difference between the phases of~$\psi_1(t,x)$ and~$\psi_2(t,x)$,
and observe that around~$t\approx0.6$, this phase difference is
indeed~$\approx\pi$.

We thus see that,
\begin{conc}
Because of the loss of phase after the collapse, the phase difference between post-collapse intersecting beams
becomes unpredictable.
\end{conc}
Therefore, as noted by Merle~\cite{Merle-92a}, the interactions between two
post-collapse beams are chaotic.

\section{Ring-type singular solutions}
 \label{sec:RingTypeSingularSolution}

In this section we
propose a continuation of ring-type singular solutions, which is
based on adding a reflecting hole with radius~$r_0$ around the
origin, and then letting~$r_0\rightarrow0+$.
\subsection{Theory
Review}\label{subsec:RingSolutionTheory} Consider the
two-dimensional, radially-symmetric critical NLS
\begin{equation}\label{eq:RadialNLS}
i\psi_t(t,r)+\psi_{rr}+\frac1{r}\psi_r
+|\psi|^2\psi=0,\qquad\psi(0,r)=\psi_0(r).
\end{equation}
Let us denote the location of the maximal amplitude by
\begin{equation}\label{eq:RmaxDefinition}
r_{max}(t)=\arg\max_r|\psi|.
\end{equation}
Singular solutions of~(\ref{eq:RadialNLS}) are called
\textit{'peak-type'} when~$r_{max}(t)\equiv0$ for~$0\leq t\leq T_c$,
and \textit{'ring-type'} when~$r_{max}(t)>0$~for~$0\leq t<T_c$.

Let
\begin{subequations}\label{eq:ExplicitSingleRingSolution}
\begin{equation}
\psi_G^{(ex)}(t,r)=\frac1{L(t)}G\left(\frac{r}{L(t)}\right)e^{i\tau(t)+i\frac{L_t}{L}\frac{r^2}{4}},
\end{equation}
where
\begin{equation}\label{eq:ExplicitSingleRingSolution_L}
L(t)=\sqrt{1-\alpha^2t},\qquad\tau=\int_0^t\frac1{L^2(s)}ds=-\frac1{\alpha^2}\ln(1-\alpha^2t),
\end{equation}
\end{subequations}
and~$G(\rho)$ is a solution of
\begin{equation}\label{eq:Qprofile_ODE}
G''(\rho)+\frac{G'}{\rho}+\left[\frac{\alpha^4}{16}\rho^2-1\right]G+G^3=0,\qquad
0\neq G(0)\in\mathbb{R},\qquad G'(0)=0.
\end{equation}
In~\cite{Gprofile-05}, Fibich, Gavish and Wang showed
that~$\psi_G^{ex}$ is an explicit ring-type solution of the
radially-symmetric critical NLS~(\ref{eq:RadialNLS}) that blows up
(in~$L^4$) at~$T_c=1/\alpha^2$. Setting~$t=0$
in~(\ref{eq:ExplicitSingleRingSolution}) gives the corresponding
initial condition
\begin{equation}\label{eq:ExplicitSingleRingSolution_IC}
\psi_G^0(r)=G(r)e^{-i(\alpha^2/8)r^2}.
\end{equation}

Equation~(\ref{eq:Qprofile_ODE}) has the two free
parameters~$\alpha$ and~$G(0)$. However, in the case of a
single-ring G profile, these two parameters are
related~\cite{SC_rings-07}. For example, in the numerical
simulations in this section, the G profile is the single-ring
solution of~(\ref{eq:Qprofile_ODE}) with
\begin{equation}\label{eq:Qprofile_ODE_IC}
G(0)\approx7.6\cdot10^{-6},\qquad\alpha\approx0.357832.
\end{equation}

\subsection{Vanishing-hole continuation}\label{subsec:VanishingHole}

The sub-threshold power continuation approach of
Section~\ref{sec:NonExplicitBeam} cannot be applied
to~$\psi_G^{(ex)}$, since these solutions have an infinite power. In
addition, this continuation cannot be applied to the ring-type
singular solutions which are in~$H^1$, since these solutions exist
only for~$P\gg \Pcr$~\cite{Gprofile-05}. Therefore, we now develop a
different continuation approach, which is based on a vanishing-hole
limit.

Let~$\psi(t,r)$ be a shrinking-ring singular solution of the
critical NLS~(\ref{eq:RadialNLS}) with an initial
condition~$\psi_0(r)$. Let us add a hole around the origin with
radius~$r_0$, and impose a Dirichlet boundary condition at~$r=r_0$,
which is equivalent to placing a reflecting conductor at~$r=r_0$. In
order for the initial condition~$\psi_G^0(r)$ to satisfy the
Dirichlet boundary condition, we slightly modify it with
a{\rev{cut-off}} function~$H_s(r/r_0)$, i.e.,\footnote{For example,
in our simulations we used the cut-off function
\[
H_s(\rho)=\left\{
\begin{array} {l l}
0 & \quad 0\leq \rho\leq 1,\\
P_5(\rho) & \quad 1<\rho<2,\\
1 & \quad \rho\geq 2,
\end{array}\right.,\quad P_5(\rho)=6\left(\rho-\frac{3}{2}\right)^5-5\left(\rho-\frac{3}{2}\right)^3+\frac{15}{8}\left(\rho-\frac{3}{2}\right)+\frac1{2}.
\]}~
\[
\psi_G^0(r)\rightarrow\psi_G^{0,r_0}(r):=\psi_G^0(r)\cdot
H_s\left(\frac{r}{r_0}\right).
\]

\begin{figure}[ht!]
\begin{center}
\scalebox{0.6}{\includegraphics{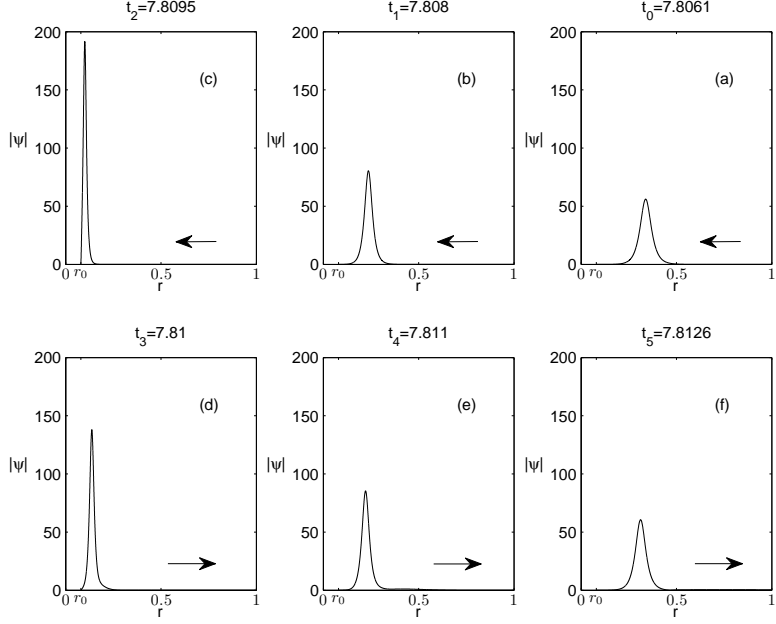}}
\caption{ Solution of the NLS~(\ref{eq:RingDirichletBCProblem})
with~$r_0=0.08$ at (a):~$t_0\approx7.806$, (b):~$t_1\approx7.808$,
(c):~$t_2\approx7.8095$ , (d):~$t_3\approx7.8100$,
(e):~$t_4\approx7.8110$, (f):~$t_5\approx 7.8126$. The arrows denote
the direction in which the ring
moves.}\label{fig:PsiQ_r0_0p08_streamDirichlet}
\end{center}
\end{figure}

We thus solve the two-dimensional, radially-symmetric critical NLS
\begin{subequations}\label{eq:RingDirichletBCProblem}
\begin{equation}
i\psi_t(t,r)+\psi_{rr}+\frac1{r}\psi_r +|\psi|^2\psi=0,\qquad
r_0<r<\infty,
\end{equation}
with the initial condition
\begin{equation}\label{eq:QprofileIC_WithHole_2}
\psi_G^{0,r_0}(r) \equiv G(r)e^{-i\frac{\alpha^2}{8}r^2}\cdot
H_s\left(\frac{r}{r_0}\right),
\end{equation}
and the Dirichlet boundary condition
\begin{equation}
\psi(t,r_0)=0,\qquad t\geq0.
\end{equation}
\end{subequations}
A typical simulation is shown in
Figure~\ref{fig:PsiQ_r0_0p08_streamDirichlet}. Initially, the ring
solution shrinks and becomes higher and narrower as it approaches
the hole. After the solution is reflected outward by the hole, it
expands and becomes lower and wider.

\begin{figure}[ht!]
\begin{center}
\scalebox{0.6}{\includegraphics{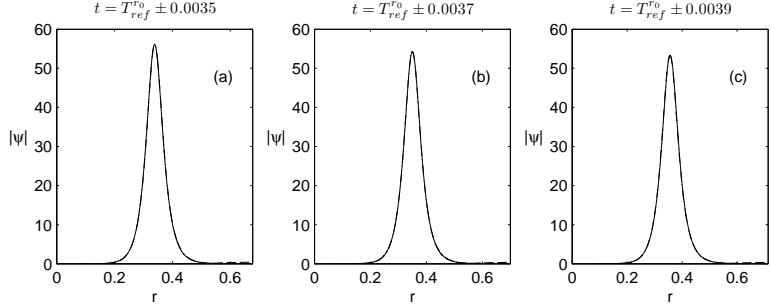}}
\caption{ Solution of Figure~\ref{fig:PsiQ_r0_0p08_streamDirichlet}
at~$t=T_{ref}^{r_0}-\Delta t$ (solid) and at~$t=T_{ref}^{r_0}+\Delta
t$ (dashes), where~$T_{ref}^{r_0}\approx7.8096$. The dotted line is
the best fitting~$\psi_G$ profile. All three curves are
indistinguishable.(a):~$\Delta t=0.0035$. (b):~$\Delta t=0.0037$.
(c):~$\Delta t=0.0039$. }\label{fig:PsiHole_GprofileFitting_2}
\end{center}
\end{figure}

We now show that the conditions of
Continuation Result~\ref{Proposition:WeakSolutionSymmetry} hold:
\begin{enumerate}
\item The solution of~(\ref{eq:RingDirichletBCProblem}) exists
globally, since otherwise it collapses at some~$0<r_c<\infty$ (a
standing ring), or at~$r\rightarrow\infty$ (an expanding ring). The
first possibility is only possible, however, if the nonlinearity is
quintic or higher, and the second possibility is not possible for
any power-nonlinearity, see~\cite{Baruch_Fibich_Gavish:2009}.
\item By continuity,~$\lim_{r_0\rightarrow0+}\psi=\psi_G^{(ex)}$ for~$0\leq
t<T_c$.
\item The solution of~(\ref{eq:RingDirichletBCProblem}) is
invariant under the transformation~(\ref{eq:NLS_SymmetryRes}).
\item $\lim_{t\rightarrow T_c}\arg\psi_G^{(ex)}=\infty$,
see~(\ref{eq:ExplicitSingleRingSolution_L}).
\item Let~$T_{ref}^{r_0}:= \arg \min_t{r_{max}(t)}$ denote the
reflection time, where~$r_{max}(t)$ is given
by~(\ref{eq:RmaxDefinition}). Since the solution focuses for~$0\leq
t<T_{ref}^{r_0}$ and defocuses for~$T_{ref}^{r_0}\leq t<\infty$, it
is collimated at~$t=T_{ref}^{r_0}$. Therefore,
\begin{equation}\label{eq:VanishingHoleCollimatedPhase}
\arg\psi(T_{ref}^{r_0},r;r_0)\equiv\alpha^{r_0},\qquad\alpha^{r_0}\in\mathbb{R}.
\end{equation}
Hence, by the arguments in the proof of
Continuation Result~\ref{Proposition:WeakSolutionSymmetry},\footnote{Here~$T_{ref}^{r_0}$
is the analog of~$T_{max}^\varepsilon$,
and~(\ref{eq:PsiRingSymmetry}) is the analog
of~(\ref{eq:PsiEpsilonSymmetry_2}).}
\begin{equation}\label{eq:PsiRingSymmetry}
\psi(T_{ref}^{r_0}+t,r;r_0)=e^{i\theta(r_0)}\psi^\ast(T_{ref}^{r_0}-t,r;r_0),\qquad
t>0.
\end{equation}
Indeed, in Figure~\ref{fig:PsiHole_GprofileFitting_2} we plot the
solution at~$t=T_{ref}^{r_0}\pm\Delta t$ for three different values
of~$\Delta t$, and observe that in all three cases, the two curves
lie on top of each other.
\end{enumerate}
Therefore, we get the following result:
\begin{prop}\label{Proposition:RingTypeWeakSol}
Let~$\psi(t,r;r_0)$ be the solution of the
NLS~(\ref{eq:RingDirichletBCProblem}), and assume
that~(\ref{eq:VanishingHoleCollimatedPhase}) holds. Then, for
any~$\theta\in\mathbb{R}$, there exists a
sequence~$r_{0,n}\rightarrow0+$ (depending on~$\theta$), such that
\begin{equation}\label{eq:AsymptoticProfileRing}
\lim_{r_{0,n}\rightarrow0+}\psi(t,r;r_{0,n})=\left\{
\begin{array}{l l}
\psi_G^{(ex)}(t,r)&\qquad0\leq t<T_c,\\
\psi^{\ast(ex)}_G(2T_c-t,r)e^{i\theta}&\qquad T_c<t,
\end{array}\right.
\end{equation}
where~$\psi_G^{(ex)}(t,r)$ is given
by~(\ref{eq:ExplicitSingleRingSolution}), and~$T_c=1/\alpha^2$.
Hence, this continuation satisfies
Properties~1 and~2.

In particular, the limiting width is given by
\[
\lim_{r_{0}\rightarrow0+}L(t;r_0)=\sqrt{\left|1-\frac{t}{T_c}\right|}.
\]
\end{prop}

\subsubsection{Simulations}
In Figure~\ref{fig:TreflectionLimit} we plot the reflection time as
a function of~$r_0$, and observe that this data is in excellent fit
with the parabola
\begin{equation}\label{eq:TrefLimit}
T_{ref}^{r_0}=\widehat{T}_c+k_2\cdot
r_0^2,\qquad\widehat{T}_c=7.80981,\qquad k_2=-0.03585.
\end{equation}
Since~$T_c=1/\alpha^2\approx7.80983$,
see~(\ref{eq:Qprofile_ODE_IC}), the extrapolation error
is~$\frac{\left|T_c-\widehat{T}_c\right|}{T_c}=0.0002\%$. The
observation that~$T_c-T_{ref}(r_0)$ scales as~$r_0^2$ and not only
as~$r_0$, will allow us to extrapolate~$\psi_G^{r_0}(t,r)$
and~$L(t;r_0)$ as a function of~$r_0^2$, rather than of~$r_0$,
leading to more accurate extrapolations.

\begin{figure}[ht!]
\begin{center}
\scalebox{0.6}{\includegraphics{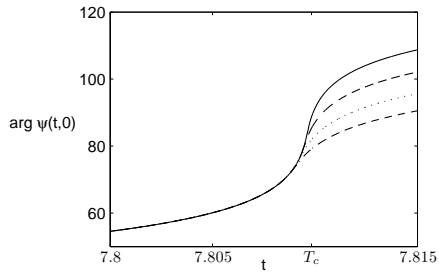}}
\caption{ Accumulated phase as a function of~$t$, for the solution
of the NLS~(\ref{eq:RingDirichletBCProblem}) with~$r_0=0.08$
(solid),~$r_0=0.1$ (dashes),~$r_0=0.125$ (dots), and~$r_0=0.15$
(dashes-dots).}\label{fig:LimitSolutionCharacterizationPhase_Dirichlet}
\end{center}
\end{figure}

We now present simulation results that support
Continuation Result~\ref{Proposition:RingTypeWeakSol}:
\begin{enumerate}
\item In Figure~\ref{fig:LimitSolutionExtrapolation_Ar}(a) we
plot~$L(t;r_0)$ for~$r_0=0.15, 0.125, 0.1$ and~$0.08$. The curve
which is obtained from the extrapolation of these curves
to~$r_0^2=0$, is nearly identical to the limiting
curve~$L=\sqrt{\left|1-\frac{t}{T_c}\right|}$, see
Figure~\ref{fig:LimitSolutionExtrapolation_Ar}(b).
\begin{figure}[ht!]
\begin{center}
\scalebox{0.6}{\includegraphics{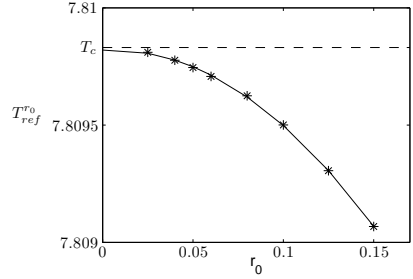}} \caption{ The
reflection time~$T_{ref}^{r_0}$ as a function of the hole
radius~$r_0$, for the solution of the
NLS~(\ref{eq:RingDirichletBCProblem}). Solid line is the
parabola~(\ref{eq:TrefLimit}).}\label{fig:TreflectionLimit}
\end{center}
\end{figure}

\begin{figure}[ht!]
\begin{center}
\scalebox{0.6}{\includegraphics{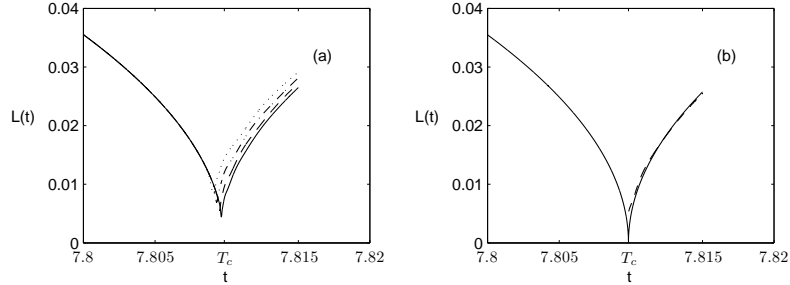}}
\caption{ Solution of the NLS~(\ref{eq:RingDirichletBCProblem}).
(a):~Solution width for~$r_0=0.15$ (dots),~$r_0=0.125$
(dashes-dots),~$r_0=0.1$ (dashes) and~$r_0=0.08$ (solid). (b):
Extrapolation of the curves~$\{L(t;r_0)\}$ from~(a) to~$r_0=0$
(dashes). Solid line is~$L(t)=\sqrt{\left|1-t/T_c\right|}$.
}\label{fig:LimitSolutionExtrapolation_Ar}
\end{center}
\end{figure}
\item In Figure~\ref{fig:LimitProfileExtrapolation_Dirichlet_Ar}(a) we
fix the time at~$t_1=7.8124>T_c$, and plot the solution profile
for~$r_0=0.15, 0.125, 0.1$ and~$0.08$. The curve which is obtained
from the extrapolation of the profiles~$\{|\psi(t_1,r;r_0)|\}$
to~$r_0^2=0$ is nearly identical to~$|\psi_G^{\ast(ex)}(2T_c-t_1,r)|$, see
Figure~\ref{fig:LimitProfileExtrapolation_Dirichlet_Ar}(b).
\begin{figure}[ht!]
\begin{center}
\scalebox{0.6}{\includegraphics{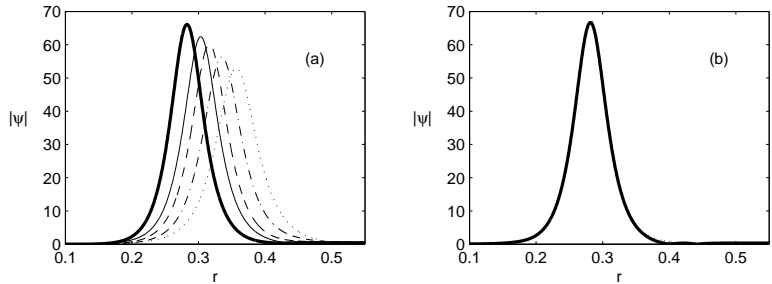}}
\caption{ (a): Solution of the NLS~(\ref{eq:RingDirichletBCProblem})
at~$t=7.8124=1.0003T_c$ for~$r_0=0.15$ (dots),~$r_0=0.125$
(dashes-dots),~$r_0=0.1$ (dashes) and~$r_0=0.08$ (solid). Solid bold
line is the extrapolation of these curves to~$r_0=0$. (b): The
extrapolated profile from (a) (solid bold). Dotted line
is~$|\psi_G^{\ast(ex)}(2T_c-7.8214,r)|$. The two curves are
indistinguishable.}\label{fig:LimitProfileExtrapolation_Dirichlet_Ar}
\end{center}
\end{figure}
\item In Figure~\ref{fig:LimitSolutionCharacterizationPhase_Dirichlet}
we plot the accumulated phase at the ring peak, i.e.,~$arg\,
\psi_G^{r_0}(t,r_{max}(t))$, and observe that small changes in~$r_0$
hardly affect the phase before the singularity, but lead
to~$\mathcal{O}(1)$ changes in the phase after the singularity.
\end{enumerate}

\section{Vanishing nonlinear-damping solutions}
\label{sec:DampedNLS}

In this section, we propose a continuation which is based on the
addition of nonlinear damping. The motivation for this approach
comes from the vanishing-viscosity solutions of hyperbolic
conservation laws. Of course, the key question is which physical
mechanism should play the role of "viscosity" in the NLS. In the
nonlinear optics context, there are numerous candidates, which
correspond to the mechanisms that are neglected in the derivation of
the NLS from Maxwell's equations: Nonparaxial effects, high-order
nonlinearities, dispersion, plasma effects, Raman, damping, etc. Of
course, for a physical mechanism to be able to play the role of
"viscosity", it should arrest the collapse regardless of how small
it is (so that we can take the limit of this term to zero, and still
have global solutions). This requirement rules out some candidates
(such as linear damping, see below), but still leaves plenty of
potential candidates (such as nonlinear saturation, see
Section~\ref{sec:nl-saturation}).

In this study we consider the case when the role of viscosity is
played by nonlinear damping. The addition of small nonlinear damping
is "physical". Indeed, in nonlinear optics, experiments suggest that
arrest of collapse is usually related to plasma formation, and
nonlinear damping can be used as a phenomenological model for the
multi-photon absorption by the plasma. In BEC, a quintic nonlinear
damping term corresponds to losses from the condensate due to
three-body inelastic recombinations. In~\cite{Markowich-04}, Bao,
Jaksch and Markowich showed numerically that the arrest of collapse
by the addition of a quintic nonlinear damping to the cubic
three-dimensional NLS is in good agreement with experimental
measurements. Nonlinear damping arises also in the context
of the complex Ginzburg-Landau equation, see Section~\ref{sec:CGL}.

\subsection{Effect of linear and nonlinear damping - review}

In~\cite{Fibich-2001},
Fibich studied asymptotically and numerically the effect of damping
on blowup in the critical NLS, and showed that when the damping is
linear, i.e.,
\begin{equation}\label{eq:CriticalLinearDampedLS}
i\psi_t(t,\textbf{x})+\Delta\psi+|\psi|^{4/d}\psi+i\delta\psi=0,\qquad\psi(0,\x)=\psi_0(\x),
\end{equation}
if the initial condition~$\psi_0(\x)$ is such that the solution
of~(\ref{eq:CriticalLinearDampedLS}) becomes singular
for~$\delta=0$, then the solution
of~(\ref{eq:CriticalLinearDampedLS}) exists globally only
if~$\delta$ is above a threshold value~$\delta_c>0$ (which depends
on~$\psi_0$). Therefore, {\em{linear damping cannot play the role of
viscosity in defining weak solutions of the NLS}}. When, however,
the damping exponent is critical or supercritical, i.e.,
\begin{equation}\label{eq:CriticalDampedNLS}
i\psi_t(t,\textbf{x})+\Delta\psi+\left(1+i\delta\right)|\psi|^{4/d}\psi=0,\qquad\delta>0,
\end{equation}
or
\begin{equation}\label{eq:SuperCriticalDampedNLS}
i\psi_t(t,\textbf{x})+\Delta\psi+|\psi|^{4/d}\psi+i\delta|\psi|^p\psi=0,\qquad\delta>0,\qquad
p>4/d,
\end{equation}
respectively, then regardless of how small~$\delta$ is, collapse is
always arrested. Therefore, Fibich suggested that nonlinear damping
can "play the role of viscosity" in defining weak NLS solutions,
i.e., we can define the continuation
\begin{equation}\label{eq:DampedNLS_Continuation}
\psi:=\lim_{\delta\rightarrow 0+}\psi^{(\delta)},
\end{equation}
where~$\psi^{(\delta)}$ is the solution
of~(\ref{eq:CriticalDampedNLS})
or~(\ref{eq:SuperCriticalDampedNLS}).

Since the results in~\cite{Fibich-2001} are not rigorous, we now
present the relevant rigorous results that exist in the literature.
{\rev{Passot, Sulem and Sulem proved that high-order nonlinear
damping always prevents collapse for~$d=2$. Antonelli and Sparber
extended this result to~$d=1$ and~$d=3$:}}
\begin{thrm}\label{thrm:SulemThrm}
~\cite{Sulem-Passot-2005,Antoneli-Sparber-2010} {\rev{The
d-dimensional cubic NLS with nonlinear damping
\begin{equation}\label{eq:CauchyProblemNLSCubic_1}
i\psi_t(t,\X)+\Delta\psi+\lambda|\psi|^2\psi+i\delta|\psi|^{p-1}\psi=0,\qquad\lambda\in\mathbb{R},\quad\delta>
0,
\end{equation}
where~$\psi_0(\x)\in H^1(\mathbb{R}^d)$,  $3<p<\infty$ if~$d=1,2$, and~$3<p<5$
if~$d=3$,
has a unique global in-time solution.}}
\end{thrm}
This rigorously shows that high-order nonlinear damping can play the
role of "viscosity". More recently, Antonelli and Sparber proved
global existence for the case where the damping exponent is equal to
that of the nonlinearity:
\begin{thrm}\label{thrm:AntonelliThrm}
~\cite{Antoneli-Sparber-2010} Consider the cubic nonlinear NLS with
a cubic nonlinear damping
\begin{equation}\label{eq:CauchyProblemNLSCubic}
i\psi_t(t,\x)+\Delta\psi+(1+i\delta)|\psi|^2\psi=0,
\end{equation}
where~$\psi_0(\x)\in H^1(\mathbb{R}^d),~\x\psi_0\in
L^2(\mathbb{R}^d)$, and~$d\leq3$. Then, for any~$\delta\geq1$,
equation~(\ref{eq:CauchyProblemNLSCubic}) has a unique global
in-time solution.
\end{thrm}
Unfortunately, because of the constraint~$\delta\geq1$,
Theorem~\ref{thrm:AntonelliThrm} does not show that critical
nonlinear damping can play the role of viscosity. We note, however,
that the asymptotic analysis and simulations of~\cite{Fibich-2001}
strongly suggest that the solution of~(\ref{eq:CriticalDampedNLS})
exists globally for any~$0<\delta\ll1$.

\subsection{Explicit continuation
when~$\psi_0(r)=\psiex(0,r)$}\label{subsec:ExplicitContDampedNLS}

In the special case where~$\psi^{(\delta)}$ is the solution
of~(\ref{eq:CriticalDampedNLS}) with the initial
condition~$\psiex(t=0)$, we can calculate explicitly the vanishing
nonlinear-damping limit~(\ref{eq:DampedNLS_Continuation}):
\begin{prop}
\label{Proposition:DampedNLSWeakSol}

Let~$\psi^{(\delta)}(t,r)$ be the solution of the
NLS~(\ref{eq:CriticalDampedNLS}) with the initial condition
\begin{equation}\label{eq:LemmaIC}
\psi_0(r)=\psiex(0,r),
\end{equation}
see~\eqref{eq:ExplicitBlowupSolution}. Then, for any~$\theta\in\mathbb{R}$, there exists a
sequence~$\delta_n\rightarrow0+$ (depending on~$\theta$), such that
\begin{equation}\label{eq:DampedWeakSolution}
\lim_{\delta_n\rightarrow 0+}\psi^{(\delta_n)}(t,r)=\left\{
\begin{array}{l l}
  \psiex(t,r) & \quad~0\leq t<T_c,\\
  \psi_{\rm explicit, \kappa}^*(2T_c-t,r)e^{i\theta} & \quad T_c<t<\infty,
\end{array} \right.
\end{equation}
where~$\psi_{\rm explicit, \kappa}$ is given
by~(\ref{eq:ExplicitBlowupSolutionAlpha_1}) with $\alpha = \kappa$,
\begin{equation}\label{eq:AlphaExpression}
\kappa=\pi\left[B_i(0)A'_i(s^*)-A_i(0)B'_i(s^*)\right]\approx 1.614,
\end{equation}
$A_i(s)$ and~$B_i(s)$ are the Airy and Bairy functions,
respectively, and~$s^*\approx-2.6663$ is the first negative root
of~$G(s)=\sqrt{3}A_i(s)-B_i(s)$.

In particular, the limiting width of the solution is given by
\begin{equation}\label{eq:DampedNLS_L_weak}
\lim_{\tilde{\delta}\rightarrow0+}L(t;\delta)=\left\{
\begin{array}{l l}
  T_c-t & \quad 0\leq t<T_c,\\
  \kappa(t-T_c) & \quad T_c<t<\infty,
\end{array} \right.
\end{equation}
\end{prop}
\Beginproof See section~\ref{sec:proof-nl-damping-psiexplicit}.
\Endproof

Equation~(\ref{eq:DampedWeakSolution}) provides a continuation of
the explicit blowup solution~$\psiex$ beyond the singularity. As
with all the continuations that we saw so far, the weak solution is
only determined up to a multiplicative phase constant~$e^{i\theta}$
(Property~1). In contrast with these
continuations, however,~{\em{the
continuation~(\ref{eq:DampedWeakSolution}) is asymmetric with
respect to~$T_c$}}, since~$\kappa\neq1$. Therefore, it does not
satisfy Property~2. This is due to the
directionality in~$t$ of the damping effect, i.e., the fact that
equation~(\ref{eq:CriticalDampedNLS}) is not invariant under the
transformation~(\ref{eq:NLSInvariantTransformations}).
\begin{rmk}
The value of $\kappa$ is also given by, see~\eqref{eq:kappa-simple},
$$
  \kappa = \frac{A_i(0)}{-A_i(s^*)}.
$$
\end{rmk}
\begin{rmk}
The limiting solution in Continuation
Result~\ref{Proposition:DampedNLSWeakSol} does not have a
non-collapsing "tail", since the power of~$\psi_{\rm explicit,
\kappa}^*(2T_c-t,r)e^{i\theta}$ is equal to that of~$\psi_0(r)$.
\end{rmk}

\begin{figure}[ht!]
\begin{center}
\scalebox{0.6}{\includegraphics{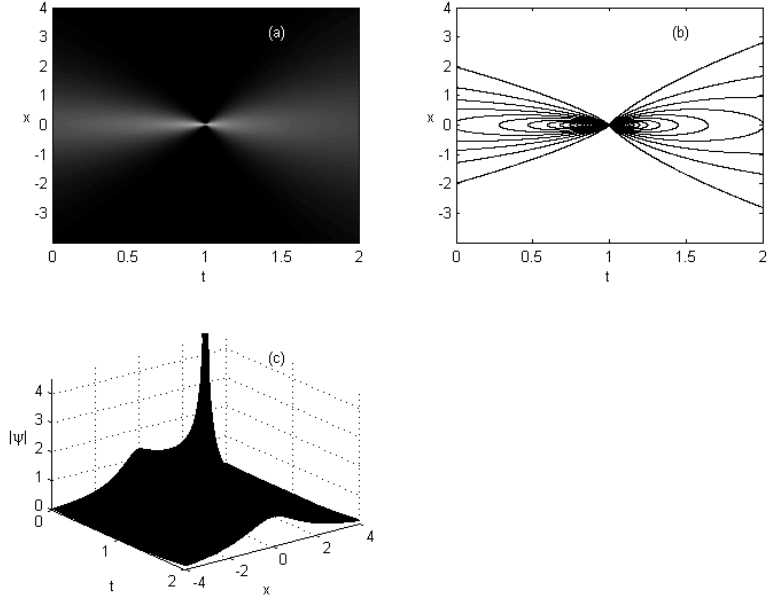}}
\caption{ The vanishing nonlinear damping
limit~(\ref{eq:DampedWeakSolution}) for~$d=1$ and~$T_c=1$. (a) Color
plot.~(b) Contour plot.~(c) Surface plot.}
\label{fig:DampedNLS_ExplicitLimitContour_new}
\end{center}
\end{figure}

\begin{rmk}\label{Remark:Lt_InfinteVelocity}
In order to understand why~$L_t^2$ increases (and not decreases)
after the singularity, we note that while the vanishing
nonlinear-damping does not affect the solution power, it increases
the Hamiltonian, see Section~\ref{sec:H-dynamics}. Since
$$
H(\psi_{\rm explicit, \alpha}) =  \frac{L_t^2}{4} ||r R^{(0)}||_2^2 =  M \alpha^2,
$$
the increase in the Hamiltonian implies that the defocusing velocity
(angle) should be higher than the focusing velocity (angle), see
Figure
~\ref{fig:DampedNLS_ExplicitLimitContour_new}.
\end{rmk}
%


\subsubsection{Simulations}

\begin{figure}[ht!]
\begin{center}
\scalebox{0.6}{\includegraphics{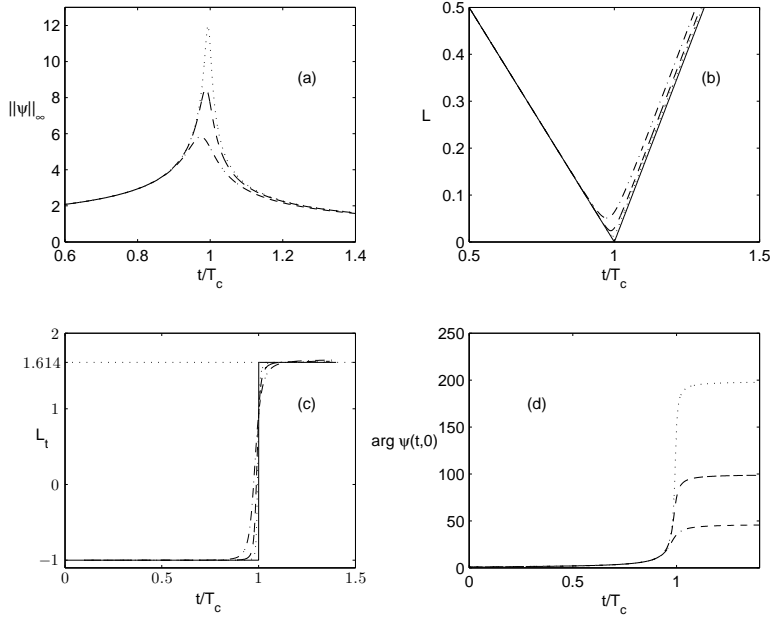}}
\caption{ Solution of the damped NLS~(\ref{eq:CriticalDampedNLS})
with~$d=1$ and the initial condition~(\ref{eq:LemmaIC}) with~$T_c=1$
for~$\delta=10^{-5}$ (dashes-dots),~$\delta=10^{-6}$ (dashes)
and~$\delta=1.25\cdot10^{-7}$ (dots). (a):~$\|\psi\|_\infty$.
(b):~$L(t;\delta)$, recovered from~$\psi$
using~(\ref{eq:BeamWidthDefinition}). Solid line
is~(\ref{eq:DampedNLS_L_weak}). (c):~$L_t(t;\delta)$. Solid line
is~(\ref{eq:DampedNLS_Lt_weak}). (d):~Accumulated phase at~$x=0$.}
\label{fig:PhaseNonUniqueness_dampedNLS}
\end{center}
\end{figure}

In order to provide numerical support to Continuation
Result~\ref{Proposition:DampedNLSWeakSol}, we solve numerically the
damped NLS~(\ref{eq:CriticalDampedNLS}) with~$d=1$ and the initial
condition~(\ref{eq:LemmaIC}) with~$T_c=1$, for various values
of~$\delta$. Figure~\ref{fig:PhaseNonUniqueness_dampedNLS}(a) shows
that as~$\delta\rightarrow0+$, the maximal amplitude increases, and
that it is attains at~$T^\delta_{max}\rightarrow T_c$.
Figures~\ref{fig:PhaseNonUniqueness_dampedNLS}(b)
and~\ref{fig:PhaseNonUniqueness_dampedNLS}(c) show
that~$\lim_{\delta\rightarrow0+}L(t;\delta)$ is given
by~(\ref{eq:DampedNLS_L_weak}), and that
\begin{equation}
   \label{eq:DampedNLS_Lt_weak}
\lim_{\tilde{\delta}\rightarrow0}L_t(t;\delta)=\left\{
\begin{array}{l l}
  -1 & \quad 0\leq t<T_c,\\
  \kappa & \quad T_c<t<\infty,
\end{array} \right.
\end{equation}
respectively, where~$\kappa$ is defined
in~(\ref{eq:AlphaExpression}).

In Figure~\ref{fig:PhaseNonUniqueness_dampedNLS}(d) we plot the
accumulated phase at~$x=0$, and observe that small changes
in~$\delta$ have a negligible effect on the phase before the
singularity, but an~$\mathcal{O}(1)$ effect on the phase after the
singularity, which is an indication that phase of the weak solution
becomes non-unique for~$t>T_c$.

\subsection{Continuation for loglog
collapse}\label{subsec:loglogCollapse}

Let~$\psi$ be a solution of the undamped critical NLS that undergoes
a loglog collapse. Since nonlinear damping leads to defocusing (and
not to oscillations) after it arrests the
collapse~\cite{Fibich-2001}, the limiting solution has a point
singularity and not a filament singularity. Indeed, the continuation
has an infinite-velocity{\rev{expanding core}}:\footnote{The
observation that the velocity of the expanding solution is infinite,
is due to Merle~\cite{Merle-private-11}.}
\begin{prop}
\label{conject:DampedNLS_Conject}
Let $\psi_0(r)$ be a radial initial condition, such that the
corresponding solution~$\psi$ of the undamped critical
NLS~(\ref{eq:dDim_CriticalNLS}) collapses with the
$\psi_{R^{(0)}}$~profile at the loglog law blowup rate at~$T_c$. Let
$\psi^{(\delta)}$ be the solution of the damped
NLS~(\ref{eq:CriticalDampedNLS}) with the same initial condition.
Then,
$$
\lim_{\delta \to 0+} \psi^{(\delta)} = \psi, \qquad 0 \le t < T_c.
$$
In addition, for any~$0<\delta\ll1$, there
exists~$\theta(\delta)\in\mathbb{R}$, and a function~$\phi\in L^2$,
such that
$$
\lim_{\delta\rightarrow
0+}\left[\psi^{(\delta)}(t,r)-{\rev{\psi^\ast_{R^{(0)}}}}(2T_c-t,r;\delta)e^{i\theta(\delta)}\right]\stackrel{L^2}\longrightarrow\phi(r),
\qquad t\longrightarrow T_c+,
$$
{\rev{where~$\psi_{R^{(0)}}$ is given
by~(\ref{eq:AsymptoticProfile}) with some function~$L(t;\delta)$,
such that}}
\[
{\rev{\lim_{t\rightarrow
T_c+}\lim_{\delta\rightarrow0+}L(t;\delta)=0,\qquad\lim_{t\rightarrow
T_c+}\lim_{\delta\rightarrow0+}L_t(t;\delta)=\infty}},
\qquad\lim_{\delta\rightarrow0+}\theta(\delta)=\infty.
\]

\end{prop}
\Beginproof See
Section~\ref{subsec:ReducedEqDampedNLS-loglog}.
\Endproof

The post-collapse infinite velocity of the expanding core, is a
consequence of the infinite velocity of the loglog collapse before
the singularity, and the increase of the velocity after the
singularity (see Remark~\ref{Remark:Lt_InfinteVelocity}).
\begin{rmk}
{\rev{Because of the infinite velocity of the expanding core, it
"immediately" interacts with the non-collapsing tail. Therefore, the
validity of the reduced equations that are used in the derivation of Continuation
Result~\ref{conject:DampedNLS_Conject} breaks down ``shortly'' after the arrest of collapse. See
Section~\ref{subsubsec:ValidOfRE} for further discussion.}}
\end{rmk}

\subsubsection{Simulations}

\begin{figure}[ht!]
\begin{center}
\scalebox{0.6}{\includegraphics{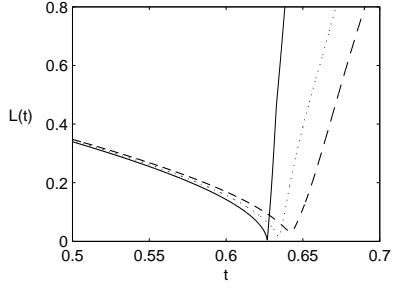}}
\caption{$L(t)$ for the solution of the damped
NLS~(\ref{eq:CriticalDampedNLS}) with~$d=1$, and the initial
condition $\psi_0(r)=\sqrt{1.05} \psiex(t=0,r;T_c=1)$
for~$\delta=10^{-3}$ (solid),~$\delta=2\cdot10^{-3}$ (dots),
and~$\delta=2.5\cdot10^{-3}$ (dashes).}
\label{fig:DampedNLSLogLog_Fit}
\end{center}
\end{figure}

\begin{figure}[ht!]
\begin{center}
\scalebox{0.6}{\includegraphics{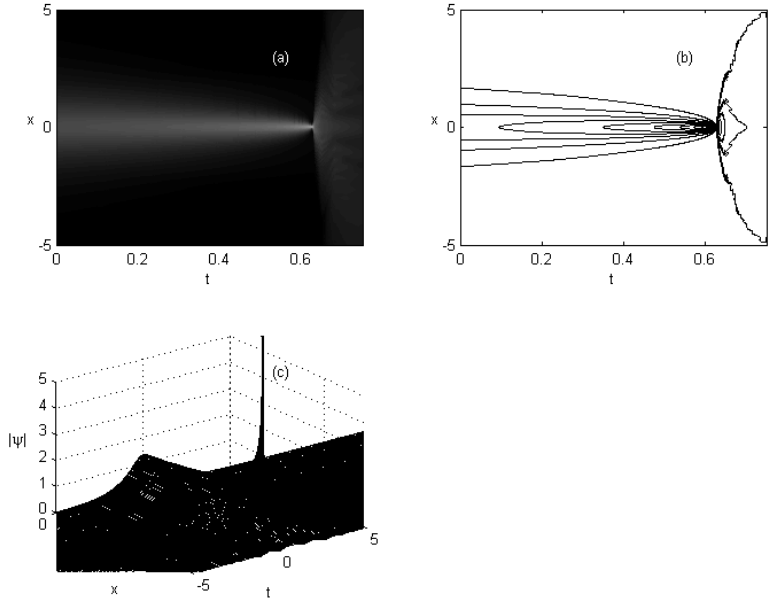}}
\caption{Same as Figure~\ref{fig:DampedNLS_ExplicitLimitContour_new}
for the numerical solution of the damped
NLS~(\ref{eq:CriticalDampedNLS}) with~$d=1$,~$\delta=10^{-3}$, and
the initial condition $\psi_0(r)=\sqrt{1.05} \psiex(t=0,r;T_c=1)$.}
\label{fig:DampedNLS_LoglogContour_new}
\end{center}
\end{figure}

In order to illustrate Continuation
Result~\ref{conject:DampedNLS_Conject} numerically,
in{\rev{Figures~\ref{fig:DampedNLSLogLog_Fit}--\ref{fig:DampedNLS_LoglogContour_new}}}
we solve the damped NLS~(\ref{eq:CriticalDampedNLS}) with~$\psi_0 =
\sqrt{1.05} \psiex(t=0)$
for~$\delta=2\cdot10^{-3}$,~$\delta=1.5\cdot10^{-3}$,
and~$\delta=10^{-3}$, and observe that the solutions are highly
asymmetric with respect to~$T_{max}^{(\delta)}$.{\rev{In addition,
as~$\delta\rightarrow0+$, the post-collapse expansion of the
singular core becomes faster and faster.}}


\subsection{Proof of Continuation Result~\ref{Proposition:DampedNLSWeakSol}}
\label{sec:proof-nl-damping-psiexplicit}

In order to prove Continuation Result~\ref{Proposition:DampedNLSWeakSol}, we
first approximate the NLS~(\ref{eq:CriticalDampedNLS}) with a
reduced system of ordinary differential equations. Then, we solve
the reduced system explicitly as~$\delta\rightarrow0+$.

\subsubsection{Reduced equations}
\label{subsec:ReducedEqDampedNLS}

\begin{lem}
Let~$\psi^{(\delta)}$ be the solution of the damped
NLS~(\ref{eq:CriticalDampedNLS}) with the initial
condition~(\ref{eq:LemmaIC}).
If~$\psi^{(\delta)}\sim\psi_{R^{(0)}}$,
see~(\ref{eq:AdiabaticCollapse}), the reduced equations for~$L(t)$
are given by
\begin{subequations}\label{eq:DampedReducedEquation}
\begin{eqnarray}
 &&\beta_t(t)=-\frac{\tilde{\delta}}{L^2},\qquad \tilde{\delta}=\frac{2c_d\delta}{M}\label{eq:DampedReducedEquation_t}, \\
 &&L_{tt}(t)=-\frac{\beta(t)}{L^3},
\end{eqnarray}
\end{subequations}
where~$c_d=\|R^{(0)}\|_{4/d+2}^{4/d+2}$, and~$M$ is given
by~(\ref{eq:ReducedEqCont}), with the initial conditions
\begin{equation}\label{eq:ReducedEqIC}
\beta(0)=0,\qquad L(0)=T_c,\qquad L_t(0)=-1.
\end{equation}
\end{lem}
\Beginproof
 In~\cite{Fibich-2001}, Fibich used {\em modulation
theory}~\cite{PNLS-99} to show that
when~$\psi_{\rm core}\sim\psi_{R^{(0)}}$, self focusing
in~(\ref{eq:CriticalDampedNLS}) is given, to leading-order, by
\begin{equation}
 \label{eq:beta_LeadingOrd}
\beta_t(t)=-\frac{\upsilon(\beta)}{L^2}-\frac{2c_d\delta}{M}\frac1{L^2},\qquad
\beta(t)=-L^3L_{tt},
\end{equation}
where~$\upsilon(\beta)$ is given by~(\ref{eq:ReducedEqCont}). We
recall that when~$\delta=0$, the initial
condition~(\ref{eq:LemmaIC}) leads to the explicit blowup
solution~(\ref{eq:ExplicitBlowupSolution}), for which~$L(t)=T_c-t$.
Therefore, when~$\delta=0$,
\begin{equation}\label{eq:L_IC}
L(0)=T_c,\qquad L_t(0)=-1,\qquad L_{tt}(0)=0.
\end{equation}
The initial condition is independent of the subsequent dynamics,
hence it is independent of~$\delta$. Therefore, the initial
condition is also given by~(\ref{eq:L_IC}) for~$\delta>0$.
Therefore, since~$\beta(t)=-L^3L_{tt}$, then~$\beta(0)=0$. Now,
since~$\upsilon(\beta=0)=0$, see~(\ref{eq:ReducedEqCont}),
by~(\ref{eq:beta_LeadingOrd}) we have that~$\beta_t<0$.
Therefore,~$\beta\leq0$, and consequently~$\upsilon(\beta)\equiv0$,
see (\ref{eq:ReducedEqCont}). \Endproof

\subsubsection{Simulations}
\label{subsec:NumericalSupportReducedEqDampedNLS}

\begin{figure}[ht!]
\begin{center}
\scalebox{0.6}{\includegraphics{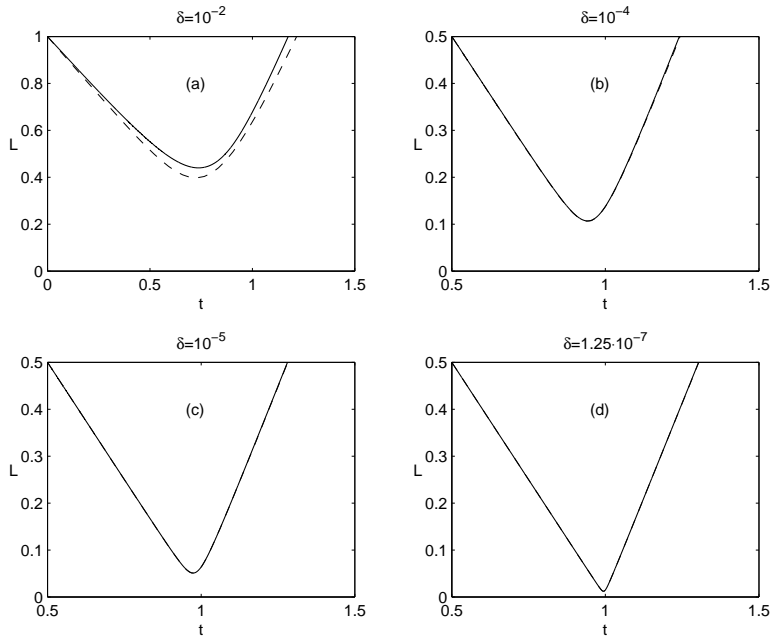}}
\caption{Width~$L(t)$ of the solution  of the damped
NLS~(\ref{eq:CriticalDampedNLS}) with~$d=1$ and the initial
condition~(\ref{eq:LemmaIC}) with~$T_c=1$ (solid). Dashed line is
the solution of the reduced
equations~(\ref{eq:DampedReducedEquation})-(\ref{eq:ReducedEqIC}).
(a):~$\delta=10^{-2}$. (b):~$\delta=10^{-4}$. (c):~$\delta=10^{-5}$.
(d):~$\delta=1.25\cdot10^{-7}$.}
\label{fig:Lexplicit_SimulationAndReducedEq}
\end{center}
\end{figure}

The derivation of the reduced
equations~(\ref{eq:DampedReducedEquation}) is based on modulation
theory, which is not rigorous. Therefore, we now provide numerical
support for the validity of~(\ref{eq:DampedReducedEquation}). In
Figure~\ref{fig:Lexplicit_SimulationAndReducedEq} we
solve~(\ref{eq:DampedReducedEquation}) numerically for~$d=1$ and
various values of~$\delta$. We compare these solutions with direct
simulations of the damped NLS~(\ref{eq:CriticalDampedNLS})
with~$d=1$ and the initial condition~(\ref{eq:LemmaIC}), from which
we extract the value of~$L(t)$ using~(\ref{eq:BeamWidthDefinition}).
When~$\delta=10^{-2}$, the two curves are similar, and
for~$\delta\leq10^{-4}$ the two curves are indistinguishable. This
confirms that as~$\delta\rightarrow 0+$, the dynamics of the
solution of the damped NLS~(\ref{eq:CriticalDampedNLS}) with the
initial condition~(\ref{eq:LemmaIC}) is given by the reduced
equations~(\ref{eq:DampedReducedEquation}).

\begin{figure}[ht!]
\begin{center}
\scalebox{0.6}{\includegraphics{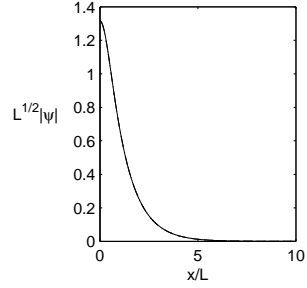}}
\caption{ Rescaled profile~$L^{1/2}(t)|\psi(t,x/L)|$ of the solution
of the damped NLS~(\ref{eq:CriticalDampedNLS}) with~$d=1$ and the
initial condition~(\ref{eq:LemmaIC}) with~$T_c=1$
and~$\delta=1.25\cdot10^{-7}$, at~$t=0.6<T_c$ (solid) and
at~$t=1.4>T_c$ (dots), where~$L(t)$ is given
by~(\ref{eq:BeamWidthDefinition}). The dashed curve is~$R(x)$. All
three curves are indistinguishable.}
\label{fig:PhaseNonUniqueness_dampedNLS_2}
\end{center}
\end{figure}

In Figure~\ref{fig:PhaseNonUniqueness_dampedNLS_2} we plot the
rescaled profile~$L^{1/2}(t)|\psi(t,x/L)|$ at~$t=0.6<T_c$ and
at~$t=1.4>T_c$. The two rescaled profiles are indistinguishable from
each other and also from~$R(x)$, providing support to the assumption
that~$\psi^{(\delta)}\sim\psi_{R^{(0)}}$, which was used in the
derivative of the reduced equations~\cite{PNLS-99}.

\subsubsection{Analysis of the reduced equations}
\label{subsec:AnalysisReducedEqDampedNLS} Our ultimate goal is to
solve the ODE system~(\ref{eq:DampedReducedEquation}) with the
initial conditions~(\ref{eq:ReducedEqIC}) explicitly
as~$\delta\rightarrow0+$. We first prove the following Lemma:
\begin{lem}\label{Lemma:AexpressionLemma}
The solution of the reduced
equations~(\ref{eq:DampedReducedEquation}) with the initial
conditions~(\ref{eq:ReducedEqIC}) can be written as
\[
L(t)=1/A(s(t)),
\]
where
\begin{equation}\label{eq:A_expression}
A(s)=\pi\left[\tilde{\delta}^{-1/3}A_i(0)-\frac1{T_c}A'_i(0)\right]\left[\sqrt{3}A_i(s)-B_i(s)\right]+\frac1{T_c}\frac{B_i(s)}{B_i(0)},
\end{equation}
$A_i(s)$ and~$B_i(s)$ are the Airy and Bairy functions,
respectively, and
\begin{equation}\label{eq:s_definition}
s(t)=-\tilde{\delta}^{1/3}\int_0^tA^2(w) \, dw.
\end{equation}
\end{lem}
\Beginproof Equation~(\ref{eq:DampedReducedEquation_t}) can be
rewritten as
\[
\beta_{\tau}=-\tilde{\delta},\qquad\tau=\int_0^t\frac1{L^2(s)}ds.
\]
Since~$\beta(t=0)=\beta(\tau=0)=0$,
\begin{equation}\label{eq:BetaEpsilonTau}
\beta(\tau)=-\tilde{\delta}\tau.
\end{equation}
We recall that~\cite{PNLS-99}
\begin{equation}\label{eq:Relation1}
\beta=\frac{A_{\tau\tau}}{A},\qquad A=\frac1{L}.
\end{equation}
Substituting~(\ref{eq:BetaEpsilonTau}) in~(\ref{eq:Relation1})
gives~$A_{\tau\tau}=-\tilde{\delta}\tau A$. The variable
change~$s=-\tilde{\delta}^{1/3}\tau$ transforms this equation into
Airy's equation
\begin{equation}\label{eq:Airy'sEquation}
A_{ss}=sA.
\end{equation}
Since~$A=1/L$, then by~(\ref{eq:ReducedEqIC})
and~(\ref{eq:s_definition}), the initial condition
for~(\ref{eq:Airy'sEquation}) are
\begin{subequations}\label{eq:IC_A_s}
\begin{eqnarray}
A(s=0)&=&A(t=0)=\frac1{T_c},\\
A_s(s=0)&=&-\tilde{\delta}^{-1/3}A_t(t=0)A^{-2}(t=0)=-\tilde{\delta}^{-1/3}.
\end{eqnarray}
\end{subequations}
The solution of Airy's equation is a linear combination of the Airy
and Bairy functions:
\begin{equation}\label{eq:Airy'sEquationSolution}
A(s)=k_1A_i(s)+k_2B_i(s).
\end{equation}
Substituting~(\ref{eq:IC_A_s}) in~(\ref{eq:Airy'sEquationSolution})
gives
\[
k_1=\frac{\tilde{\delta}^{-1/3}B_i(0)+\frac1{T_c}B'_i(0)}{A_i(0)B'_i(0)-A'_i(0)B_i(0)},
\quad
k_2=-\frac{\tilde{\delta}^{-1/3}A_i(0)+\frac1{T_c}A'_i(0)}{A_i(0)B'_i(0)-A'_i(0)B_i(0)}.
\]
The Airy and Bairy functions satisfy the Wronskian relation
\begin{equation}\label{eq:AiryWronskianRelation}
A_i(s)B'_i(s)-A'_i(s)B_i(s)=\frac1{\pi}.
\end{equation}
Therefore,
\begin{equation}\label{eq:Airy'sConstants}
k_1=\pi\left(\tilde{\delta}^{-1/3}B_i(0)+\frac1{T_c}B'_i(0)\right),
\quad
k_2=-\pi\left(\tilde{\delta}^{-1/3}A_i(0)+\frac1{T_c}A'_i(0)\right).
\end{equation}
The Airy and Bairy function satisfy
\begin{equation}\label{eq:AiryBairyRelations}
A_i(0)=B_i(0)/\sqrt{3},\qquad A'_i(0)=-B'_i(0)/\sqrt{3}.
\end{equation}
Substituting of~(\ref{eq:Airy'sConstants})
in~(\ref{eq:Airy'sEquationSolution}), using
relations~(\ref{eq:AiryBairyRelations})
and~(\ref{eq:AiryWronskianRelation}) leads
to~(\ref{eq:A_expression}). \Endproof
\subsubsection{$\tilde{\delta}\rightarrow0+$}
\begin{lem}\label{Lemma:t_s_limit}
Let~$s^*_{\tilde{\delta}}$ be the first negative root
of~$A(s;\tilde{\delta})$, see~(\ref{eq:A_expression}). Then,
$\lim_{s\rightarrow s^*_{\tilde{\delta}}}t(s;\delta)=\infty$.
\end{lem}
\Beginproof Inverting~(\ref{eq:s_definition}) gives
\[
t(s^*_{\tilde{\delta}})=\tilde{\delta}^{-1/3}\int^0_{s^*_{\tilde{\delta}}}\frac1{A^2(s)}ds.
\]
Let~$0<\epsilon\ll1$. Then,
\[
t(s^*_{\tilde{\delta}})>\tilde{\delta}^{-1/3}\int^{s^*_{\tilde{\delta}}+\epsilon}_{s^*_{\tilde{\delta}}}\frac1{A^2(s)}ds.
\]
By definition,~$A(s^*_{\tilde{\delta}};\tilde{\delta})=0$.
Furthermore, since~$A(s;\tilde{\delta})$ is a nontrivial solution of
Airy's equation,
then~$A_s(s^*_{\tilde{\delta}};\tilde{\delta})\neq0$, since
otherwise from uniqueness it follows
that~$A(s;\tilde{\delta})\equiv0$. Therefore, there
exists~$0<\varepsilon$ such that
\begin{equation}\label{eq:A_s_taylor}
A(s;\tilde{\delta})\sim(s-s^*_{\tilde{\delta}})A_s(s^*_{\tilde{\delta}};\tilde{\delta}),\qquad
s^*_{\tilde{\delta}}\leq s\leq s^{*}_{\tilde{\delta}}+\varepsilon.
\end{equation}
Hence,
\[
\int^{s^*_{\tilde{\delta}}+\epsilon}_{s^*_{\tilde{\delta}}}\frac1{A^2(s)}ds\sim
\frac1{A^2_s(s^*_{\tilde{\delta}};\tilde{\delta})}\int^{s^*_{\tilde{\delta}}+\epsilon}_{s^*_{\tilde{\delta}}}\frac1{(s-s^*_{\tilde{\delta}})^2}ds=\infty.
\]
Therefore,~$t(s^*_{\tilde{\delta}})=\infty$. \Endproof

\begin{figure}[ht!]
\begin{center}
\scalebox{0.6}{\includegraphics{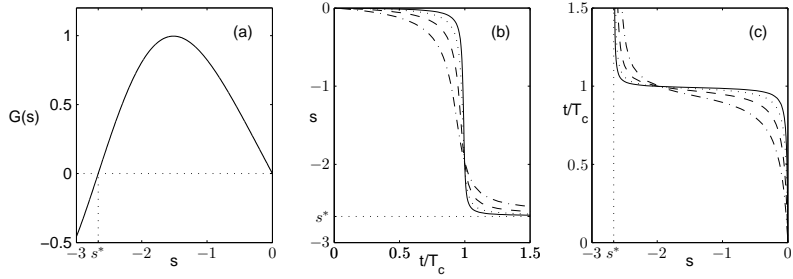}}
\caption{(a): The function~$G(s)$, see~(\ref{eq:G_s}). Here,~$s^*$
is the first negative root of~$G(s)$. (b):~$s(t)$, as calculated
from numerical integration of~(\ref{eq:s_definition})
for~$\delta=10^{-4}$ (dashes-dots),~$\delta=10^{-5}$
(dashes),~$\delta=10^{-6}$ (dots) and~$\delta=1.25\cdot10^{-5}$
(solid). (c): Same as (b) for the inverse
function~$t(s)$.}\label{fig:SstarObservations_new}
\end{center}
\end{figure}

By~(\ref{eq:s_definition}),~$s(t)$ is monotonically decreasing
from~$s(t=0)=0$. Hence, the interval~$0\leq t<\infty$ transforms
to~$0\geq s>s^*_{\delta}$. Since the Airy and the Bairy functions
are bounded for~$s\leq0$, then~$A(s;\tilde{\delta})$ is finite
for~$s\leq0$, see~(\ref{eq:A_expression}). This shows that the
solution of the damped NLS~(\ref{eq:CriticalDampedNLS}) does not
collapse. Note that~$A(s=s^*_{\tilde{\delta}})=0$ corresponds to an
infinite beam width, i.e., to a complete diffraction.
\begin{lem}\label{Lemma:Sstar_limit}
Let~$s^*=\lim_{\tilde{\delta}\rightarrow0+}s^*_{\tilde{\delta}}$.
Then,~$s^*\approx-2.6663$ is the first negative root of
\begin{equation}\label{eq:G_s}
G(s)=0,\qquad G(s):=\sqrt{3}A_i(s)-B_i(s).
\end{equation}
\end{lem}
\Beginproof By~(\ref{eq:A_expression}),
as~$\tilde{\delta}\rightarrow0+$,
\begin{equation}\label{eq:A_Sstar}
A(s)\sim\pi
A_i(0)\tilde{\delta}^{-1/3}\left[\sqrt{3}A_i(s)-B_i(s)\right] = \pi
A_i(0)\tilde{\delta}^{-1/3} G(s).
\end{equation}
Therefore,~$s^*$ satisfies~(\ref{eq:G_s}). The value of~$s^*$ was
computed numerically, see Figure~\ref{fig:SstarObservations_new}(a).
\Endproof

In Figure~\ref{fig:SstarObservations_new}(b) we plot~$s(t)$ by
numerically integrating~(\ref{eq:s_definition}), where A(s) is given
by~(\ref{eq:A_expression}), and observe that
as~$\tilde{\delta}\rightarrow0+$, ~$s(t)$ tends to the step function
\[
\lim_{\tilde{\delta}\rightarrow0+}s(t)=\left\{
\begin{array}{l l}
  0 & \quad 0\leq t<T_c,\\
  s^* & \quad T_c<t<\infty.
\end{array} \right.
\]
Therefore, the inverse function~$t(s)$ tends to the step function,
\[
\lim_{\tilde{\delta}\rightarrow0+}t(s)=\left\{
\begin{array}{l l}
\infty & \quad s=s^*,\\
  T_c & \quad s^*<s<0,\\
  0 & \quad s=0,
\end{array} \right.
\]
see Figure~\ref{fig:SstarObservations_new}(c). We now prove these
limits analytically.
\begin{lem}\label{Lemma:T0_limit}
For any~$s_1$ such that~$-1\ll
s_1<0$,~$\lim_{\tilde{\delta}\rightarrow0+}t(s_1)=T_c$.
\end{lem}
\Beginproof Equation~(\ref{eq:A_expression}) can be rewritten as
\[
A(s;\tilde{\delta})=R(s)\tilde{\delta}^{-1/3}+\frac1{T_c}F(s),
\]
where
\begin{equation}\label{eq:R_F_definition}
R(s)=\pi A_i(0)G(s),\qquad F(s)=\frac{B_i(s)}{B_i(0)}-\pi
A'_i(0)G(s).
\end{equation}
Therefore,
\begin{equation}\label{eq:integ_dt_ds}
t(s_1)=\tilde{\delta}^{-1/3}\int_{s_1}^{0}\frac1{A^2(s)}ds=\int_{s_1}^{0}\frac{\tilde{\delta}^{1/3}}{\left[R(s)+\frac{\tilde{\delta}^{1/3}}{T_c}F(s)\right]^2}ds.
\end{equation}
By~(\ref{eq:AiryBairyRelations}),~(\ref{eq:G_s})
and~(\ref{eq:R_F_definition}),
\begin{equation}\label{eq:R_F_approx_1}
R(0)=0,\qquad F(0)=1,
\end{equation}
and
\begin{equation}\label{eq:R0_tag}
R'(0)=\pi A_i(0)\left(\sqrt{3}A'_i(0)-B'_i(0)\right)=-2\pi
A_i(0)B'_i(0).
\end{equation}
Furthermore, by~(\ref{eq:AiryBairyRelations}),
\begin{subequations}
\begin{equation}\label{eq:Rel_1}
A_i(0)B'_i(0)+A'_i(0)B_i(0)=0,
\end{equation}
and by the Wronskian relation~(\ref{eq:AiryWronskianRelation}),
\begin{equation}\label{eq:Rel_2}
A_i(0)B'_i(0)-A'_i(0)B_i(0)=\frac1{\pi}.
\end{equation}
\end{subequations}
Adding~(\ref{eq:Rel_1}) to~(\ref{eq:Rel_2}), gives~$2\pi
A_i(0)B'_i(0)=1$. Therefore, by~(\ref{eq:R0_tag}),
\begin{equation}
 \label{eq:R_F_approx_2}
R'(0)=-1.
\end{equation}
Thus, by~(\ref{eq:R_F_approx_1}) and~(\ref{eq:R_F_approx_2}),
\[
R(s)\sim -s,\qquad F(s)\sim 1,\qquad-1\ll s\leq0.
\]
Substituting this in~(\ref{eq:integ_dt_ds}) gives
\begin{equation}\label{eq:t_s1}
t(s_1)\sim\int_{s_1}^{0}\frac{\tilde{\delta}^{1/3}}{\left[-s+\frac{\tilde{\delta}^{1/3}}{T_c}\right]^2}ds=\left[T_c-\frac{\tilde{\delta}^{1/3}}{-s_1
+\frac{\tilde{\delta}^{1/3}}{T_c}}\right].
\end{equation}
Therefore,~$\lim_{\tilde{\delta}\rightarrow 0+}t(s_1)=T_c$.
\Endproof

From~(\ref{eq:t_s1}) it follows that
\[
\lim_{\tilde{\delta}\rightarrow0+}t(s=-c\tilde{\delta}^{1/3})=T_c-\frac1{c+\frac1{T_c}}=\frac{T_c}{1+\frac1{c}\frac1{T_c}}.
\]
Therefore,~$t(s)$ has a boundary layer at~$s=0$ with
thickness~$\tilde{\delta}^{1/3}$, in which it increases
monotonically from 0 to~$T_c$. Hence,
\begin{corol}\label{Corollary:t_s_BoundaryLayer}
$\lim_{\tilde{\delta}\rightarrow0+}s(t;\tilde{\delta})=0,\qquad0\leq
t<T_c$.
\end{corol}
\begin{lem}\label{Lemma:limit_dt_delta}
For
any~$s^*<s_2<s_1<0$,~$\lim_{\tilde{\delta}\rightarrow0}\left(t(s_2)-t(s_1)\right)=0$.
\end{lem}
\Beginproof By~(\ref{eq:integ_dt_ds}),
\[
t(s_2)-t(s_1)=\int_{s_2}^{s_1}\frac{\tilde{\delta}^{1/3}}{\left[R(s)+\frac{\tilde{\delta}^{1/3}}{T_c}F(s)\right]^2}ds.
\]
Since~$R(s)=\pi A_i(0)G(s)$ and~$G(s)>0$ in~$(s^*,0)$,
then~$R(s)\geq c>0$ for any~$s^*<s_2\leq s\leq s_1<0$, where~$c$ is
independent of~$\tilde{\delta}$. Therefore,
as~$\tilde{\delta}\rightarrow0+$,
\[
\lim_{\tilde{\delta}\rightarrow
0}(t(s_2)-t(s_1))\sim\lim_{\tilde{\delta}\rightarrow
0}\tilde{\delta}^{1/3}\int_{s_2}^{s_1}\frac1{R^2(s)}ds=0.
\]
\Endproof

From Lemma~\ref{Lemma:T0_limit} and Lemma~\ref{Lemma:limit_dt_delta}
it follows that
\begin{corol}\label{Corollary:TransitionTime}
\[
\lim_{\tilde{\delta}\rightarrow0}t(s)=T_c,\qquad s^*<s<0.
\]
\end{corol}
\begin{lem}\label{Lemma:s_t_limit_2}
$\lim_{\tilde{\delta}\rightarrow0+}s(t;\tilde{\delta})=s^*$
for~$T_c<t<\infty$.
\end{lem}
\Beginproof Since~$s(t;\tilde{\delta})$ is monotonically decreasing
(see equation~(\ref{eq:s_definition})),~$t(s;\tilde{\delta})\approx
T_c$ for~$s^*_{\tilde{\delta}}<s<0$
(Corollary~\ref{Corollary:TransitionTime}), and~$\lim_{s\rightarrow
s^*_{\tilde{\delta}}}t(s;\tilde{\delta})=\infty$
(Lemma~\ref{Lemma:t_s_limit}), then near~$s^*_{\tilde{\delta}}$
there is a boundary layer in which~$t$ changes from~$T_c$
to~$\infty$. Therefore, the values~$(T_c,\infty)$ are attained in
the boundary layer around~$s^*_{\tilde{\delta}}$.
Since~$\lim_{\tilde{\delta}\rightarrow0+}s^*_{\tilde{\delta}}=s^*$,
the result follows. \Endproof

From Corollary~\ref{Corollary:t_s_BoundaryLayer} and
Lemma~\ref{Lemma:s_t_limit_2} it follows that, see
Figure~\ref{fig:SstarObservations_new}(b),
\begin{corol}\label{Corollary:s_t_limit}
\begin{equation}\label{eq:s_limit}
\lim_{\tilde{\delta}\rightarrow 0+}s(t;{\tilde{\delta}})=\left\{
\begin{array}{l l}
  0 & \quad 0\leq t<T_c,\\
  s^* & \quad T_c<t<\infty.
\end{array} \right.
\end{equation}
\end{corol}
Our next observation concerns~$L(t)$.
\begin{lem}\label{Lemma:L_limit}
\begin{equation}\label{eq:L_limit_DampedNLS}
\lim_{\tilde{\delta}\rightarrow0+}L(t)=\left\{
\begin{array}{l l}
  T_c-t & \quad 0\leq t<T_c,\\
  \kappa(t-T_c) & \quad T_c<t<\infty,
\end{array} \right.
\end{equation}
where~$\kappa$ is defined in equation~(\ref{eq:AlphaExpression}).
\end{lem}
\Beginproof Using~$L=A^{-1}$ and~(\ref{eq:Airy'sEquationSolution}),
\begin{equation}\label{eq:Lt_ByAiryEq}
L_t=-A^{-2}A_s\frac{ds}{dt}=-A^{-2}A_s(-\tilde{\delta}^{1/3}A^2)=\tilde{\delta}^{1/3}A_s=\tilde{\delta}^{1/3}\left[k_1A'_i(s)+k_2B'_i(s)\right].
\end{equation}
By~(\ref{eq:Airy'sConstants}),
\[
\lim_{\tilde{\delta}\rightarrow0+}\tilde{\delta}^{1/3}k_1=\pi
B_i(0),\qquad
\lim_{\tilde{\delta}\rightarrow0+}\tilde{\delta}^{1/3}k_2=-\pi
A_i(0).
\]
Therefore,
\[
\lim_{\tilde{\delta}\rightarrow0}L_t(t)=\pi\left[B_i(0)A'_i(\tilde{s}(t))-A_i(0)B'_i(\tilde{s}(t))\right],\qquad
\tilde{s}(t)=\lim_{\tilde{\delta}\rightarrow0}s(t;\tilde{\delta}).
\]
Hence, by~(\ref{eq:s_limit}),
\begin{equation*}\label{eq:Lt_limit}
\lim_{\tilde{\delta}\rightarrow0}L_t(t)=\left\{
\begin{array}{l l}
  \pi\left[B_i(0)A'_i(0)-A_i(0)B'_i(0)\right] & \quad 0\leq t<T_c,\\
  \pi\left[B_i(0)A'_i(s^*)-A_i(0)B'_i(s^*)\right] & \quad
  T_c<t<\infty.
\end{array} \right.
\end{equation*}
Since~$\pi\left[B_i(0)A'_i(0)-A_i(0)B'_i(0)\right]=-1$, see
equation~(\ref{eq:AiryWronskianRelation}),
\begin{equation*}\label{eq:Lt_limit_calculated}
\lim_{\tilde{\delta}\rightarrow0}L_t(t)=\left\{
\begin{array}{l l}
  -1 & \quad 0\leq t<T_c,\\
  \kappa & \quad T_c<t<\infty,
\end{array} \right.
\end{equation*}
where~$\kappa$ is given by~(\ref{eq:AlphaExpression}).
Since~$L(0)=T_c$,~(\ref{eq:L_limit_DampedNLS}) follows. \Endproof

Our last observation concerns the solution phase. By
definition~(\ref{eq:AsymptoticProfile}),
\[
\tau(t)=\int_0^t\frac1{L^2(s)}ds.
\]
Therefore, by~(\ref{eq:L_limit_DampedNLS}),
\begin{equation}\label{eq:TauLimit_DampedNLS}
\lim_{\tilde{\delta}\rightarrow0+}\tau(T_c)=\int_0^{T_c}\frac1{(T_c-s)^2}ds=\infty.
\end{equation}
Therefore, the phase information is lost at the singularity.

\subsubsection{Proof of
Continuation Result~\ref{Proposition:DampedNLSWeakSol}}

We have that~$\psi^{\delta}\sim\psi_{R^{(0)}}$, where~$L(t)$ is
given by~(\ref{eq:L_limit_DampedNLS}). Therefore, by
Lemma~\ref{Lemma:AsymptoticProfile_L}, when~$0\leq
t<T_c$,~$\psi_{R^{(0)}}(r,t)=\psiex(t,r)$, and
when~$T_c<t$,~$\psi_{R^{(0)}}(t,r)=\psiexalpha^\ast(2T_c-t,r)$. In
addition, since~$\arg\psi(t,0)\sim\arg\psi_{R^{(0)}}(t,0)=\tau(t)$,
equation~(\ref{eq:TauLimit_DampedNLS}) shows that the limiting phase
becomes infinite at and after the singularity. Hence, for a given
given~$t>T_c$ and~$\theta\in\mathbb{R}$, there exists a
sequence~$\delta_n\rightarrow0+$, such
that~$\lim_{\delta_n\rightarrow0+}\arg\psi^{(\delta_n)}(t,0)=\theta$.
Hence, Continuation Result~\ref{Proposition:DampedNLSWeakSol} follows.

\subsection{Hamiltonian dynamics}
   \label{sec:H-dynamics}

In the case of a non-conservative perturbation such as nonlinear damping,
the Hamiltonian of~$\psi$ can be approximated with, see~\cite[eq.~(H.5)]{PNLS-99},
\begin{equation}
\label{eq:H(psi)}
 H(\psi) \sim \frac{M}{2} (L^2)_{tt}.
\end{equation}
Since $\beta_t = -\frac{1}{2} L^2 (L^2)_{ttt}$, then
$$
H_t \sim \frac{M}{2} (L^2)_{ttt} = -M \frac{\beta_t}{L^2} =
M \frac{\tilde{\delta}}{L^4}.
$$
Therefore, the Hamiltonian {\em increases} with~$t$.
Moreover, by~\eqref{eq:s_definition} and~\eqref{eq:A_Sstar},
$$
H_s = H_t \frac{dt}{ds} =M\frac{\tilde{\delta}}{L^4} \frac{1}{-\tilde\delta^{1/3} A^2}
=- M\tilde{\delta}^{2/3} A^2 \sim -M\pi^2 A_{i}^2(0) G^2(s).
$$
Therefore,
$$
\Delta H:= H(t=\infty)-H(0) = \int_{s=0}^{s^*}H_s \, ds = M \pi^2 A_{i}^2(0)
\int_{s=s^*}^0  G^2(s)\, ds.
$$
Now,
$$
\int_{s=s^*}^0  G^2(s)\, ds = G^2 s\Big|_{s=s^*}^0 -  \int_{s=s^*}^0 2 s G G_s(s)\, ds
=- 2 \int_{s=s^*}^0 s G G_s(s)\, ds.
$$
Since $G_{ss} = sG$, then
$$
\int_{s=s^*}^0  G^2(s)\, ds = - 2 \int_{s=s^*}^0 G_{ss} G_s(s)\, ds  = -G_s^2\Big|_{s=s^*}^0.
$$
Therefore,
$$
\Delta H =-M \pi^2 A_{i}^2(0) G_s^2\Big|_{s=s^*}^0.
$$

Since the Wronskian of the Airy equation is a constant,
$$
W(G,A) = G(s) A_{i}'(s)-A_{i}(s) G'(s) \equiv  G(0) A_{i}'(0)-A_{i}(0) G'(0).
$$
Therefore, since $G$ vanishes at $s=0,s^*$,
$$
G'(s^*) = \frac{A_{i}(0) G'(0)}{A_{i}(s^*)}.
$$
Also, by~\eqref{eq:R_F_approx_2},
$$
-1 = R'(0) = \pi A_{i}(0) G'(0).
$$
Therefore,
$$
\Delta H = -M \pi^2 A_{i}^2(0) G_s^2(0)  \left( 1-\left(\frac{A_i(0)}{A_i(s^*)} \right)^2 \right) = M\left(\left(\frac{A_i(0)}{A_i(s^*)} \right)^2 - 1 \right).
$$

On the other hand, by~\eqref{eq:H(psi)}, $H \sim M(L_t^2-\beta/L^2)$. In addition, since the limiting solution has exactly the critical power, then
$\lim_{\tilde\delta \to 0+} \beta = 0$, see~\eqref{eq:BetaDefGaussian}. Therefore, $\lim_{\tilde\delta \to 0+} H  =  M L_t^2$.
Hence, by~\eqref{eq:DampedNLS_Lt_weak},
$$
\Delta H = M(\kappa^2-1).
$$
Comparison of the last two expressions for~$\Delta H$ shows that
\begin{equation}
   \label{eq:kappa-simple}
  \kappa = \frac{A_i(0)}{|A_i(s^*)|}.
\end{equation}

\subsection{Proof of Continuation Result~\ref{conject:DampedNLS_Conject}}
\label{subsec:ReducedEqDampedNLS-loglog}
\subsubsection{Analysis of the reduced equations}
{\rev{We first analyze the dynamics of the collapsing core under the
assumption that it is governed by the reduced
equations~\eqref{eq:beta_LeadingOrd}}}.\footnote{{\rev{The validity
of this assumption is discussed in
Section~\ref{subsubsec:ValidOfRE}.}}} By continuity, as $\delta
\longrightarrow 0+$, the limiting solution undergoes a loglog
collapse as $t \longrightarrow T_c-$. Therefore, as $t
\longrightarrow T_c-$, the amount of power that collapses into the
singularity is exactly~$\Pcr$. Hence, by~\eqref{eq:BetaDefGaussian},
$\lim_{\delta \to 0+} \beta(T_c) = 0$. Therefore, since $\beta_t<0$,
see~\eqref{eq:beta_LeadingOrd}, and since $\nu(\beta) \equiv 0$ for
$\beta<0$, we have that
 $\lim_{\delta \to 0+} \beta_t = 0$ for $t>T_c$.  Therefore, $\lim_{\delta \to 0+} \beta \equiv 0-$ for $t>T_c$.
Hence, the collapsing core expands linearly with the
$\psi_{R^{(0)}}$~profile. Therefore, the expanding core is given
by~$\psiexalpha^\ast(2T_c-t,r)$.

The above arguments show that{\rev{according to the reduced
equations}}, after the singularity the solution is of a
Bourgain-Wang type, but do not provide the value of the expansion
velocity~$\alpha$. In order to do so, we now solve the reduced
equations. By~(\ref{eq:beta_LeadingOrd}),
\begin{equation}\label{eq:beta_tau_Loglog}
\beta_\tau=-\nu(\beta)-\tilde{\delta},\qquad\tau=\int_0^t\frac1{L^2(s)}ds.
\end{equation}
For solutions that undergoes a loglog collapse we have
that~$\beta(0)>0$, see~(\ref{eq:BetaDefGaussian}). For a
fixed~$\tilde{\delta}>0$, as~$\beta\searrow0+$,~$\nu(\beta)$ becomes
negligible compared with~$\tilde{\delta}$. Therefore, to leading
order,~$\beta_\tau=-\tilde{\delta}$.\footnote{{\rev{In other words,
we approximate~(\ref{eq:beta_LeadingOrd})
with~(\ref{eq:DampedReducedEquation}). The validity of
neglecting~$\nu(\beta)$ is discussed in
Section~\ref{subsubsec:ValidOfRE}}}.} Hence,
\begin{equation}\label{eq:BetaTau}
\beta(\tau)=\beta_0-\tilde{\delta}\tau,\qquad\beta_0=\beta(0)>0.
\end{equation}

Substituting~(\ref{eq:BetaTau}) in~(\ref{eq:Relation1})
gives~$A_{\tau\tau}=(\beta_0-\tilde{\delta}\tau)A$. The variable
change
\begin{equation}\label{eq:s_def}
s=\tilde{\delta}^{-2/3}\beta_0-\tilde{\delta^{1/3}}\tau
\end{equation}
transforms this equation into Airy's
equation~(\ref{eq:Airy'sEquation}).

Let~$A(t=0)=A_0=1/L_0$ and~$A_t(t=0)=A_0'=-L_0'/(L_0)^2$ be the
initial conditions for~$A(t)$. Therefore, the initial conditions
for~$A(s;\tilde{\delta})$ are
\begin{subequations}\label{eq:Airy'sEqIC}
\begin{eqnarray}\label{eq:Airy'sEqIC_1}
A(s=s_0)&=&A_0,\\\label{eq:Airy'sEqIC_1}
A_s(s=s_0)&=&-\tilde{\delta}^{-1/3}A_t(t=0)A^{-2}(t=0)=-\tilde{\delta}^{-1/3}\frac{A_0'}{(A_0)^2},\label{eq:Airy'sEqIC_2}
\end{eqnarray}
where
\begin{equation}\label{eq:s0_def}
s_0:=\tilde{\delta}^{-2/3}\beta_0.
\end{equation}
\end{subequations}
Therefore,~$A(s;\tilde{\delta})$ is given
by~(\ref{eq:Airy'sEquationSolution}) with
\begin{equation}\label{eq:k1_k2}
k_1=\pi\left(\frac{A_0'}{(A_0)^2}\tilde{\delta}^{-1/3}B_i(s_0)+A_0B_i'(s_0)\right),\quad
k_2=-\pi\left(\frac{A_0'}{(A_0)^2}\tilde{\delta}^{-1/3}A_i(s_0)+A_0A_i'(s_0)\right).
\end{equation}

\subsubsection{$\tilde{\delta}\rightarrow0+$} 
\begin{figure}[ht!]
\begin{center}
\scalebox{0.6}{\includegraphics{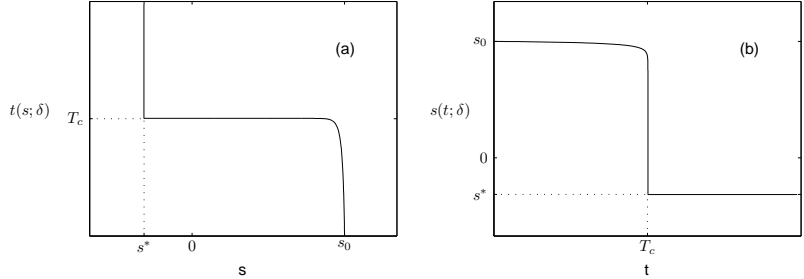}}
\caption{ (a)~$t(s;\delta)$.~(b)~$s(t;\delta)$. For
~$0<\tilde{\delta}\ll1$.}
\label{fig:t_s__s_t_DampedNLSLoglogLawWithoutNuBeta}
\end{center}
\end{figure}

In
Figure~\ref{fig:t_s__s_t_DampedNLSLoglogLawWithoutNuBeta}(a) we
plot~$t(s;\tilde{\delta})$, and observe that
\begin{lem}\label{Lemma:t_s_Lemma}
\begin{equation}\label{eq:t_s_limit_Loglog}
\lim_{\tilde{\delta}\rightarrow0+}
t\left(s;\tilde{\delta}\right)=\left\{
\begin{array}{l l}
  0 & s=s_0,\\
  T_c & \quad s^\ast<s~\mbox{and}~s_0-s\gg1\\
  \infty & \quad s=s^\ast.\\
\end{array} \right.
\end{equation}
\end{lem}
\Beginproof By~(\ref{eq:beta_tau_Loglog})
and~(\ref{eq:s_def}),~$\frac{ds}{dt}=-\tilde{\delta}^{1/3}A^2$.
Therefore,
\begin{equation}\label{eq:t_s_def_1}
t\left(s;\tilde{\delta}\right)=-\frac1{\tilde{\delta}^{1/3}}\int_{s_0}^s\frac1{A^2(s)}ds=-\frac1{\tilde{\delta}^{1/3}}\int_{s_0}^s\frac1{[k_1A_i(s)+k_2B_i(s)]^2}ds.
\end{equation}
Following~\cite{Albright-Gavathas-86}, it is easy to verify by
direct differentiation and using the Wronskian
relation~(\ref{eq:AiryWronskianRelation}) that
\begin{equation}\label{eq:t_s_def_2}
t\left(s;\tilde{\delta}\right)=\frac{\pi}{k_1\tilde{\delta}^{1/3}}\left[\frac{B_i(s_0)}{A(s_0;\tilde{\delta})}-\frac{B_i(s)}{A(s;\tilde{\delta})}\right]
=-\frac{\pi}{k_2\tilde{\delta}^{1/3}}\left[\frac{A_i(s_0)}{A(s_0;\tilde{\delta})}-\frac{A_i(s)}{A(s;\tilde{\delta})}\right].
\end{equation}
Therefore, by~(\ref{eq:Airy'sEqIC_1}),
\begin{equation}\label{eq:t_s_def_2}
t\left(s;\tilde{\delta}\right)=\frac{\pi}{k_1\tilde{\delta}^{1/3}}\left[\frac{B_i(s_0)}{A_0}-\frac{B_i(s)}{A(s;\tilde{\delta})}\right]
=\frac{\pi}{k_2\tilde{\delta}^{1/3}}\left[\frac{A_i(s)}{A(s;\tilde{\delta})}-\frac{A_i(s_0)}{A_0}\right].
\end{equation}
\begin{lem}\label{Lemma:t_s_limit_}
Let~$s^*_{\tilde{\delta}}$ be the first negative root
of~$A(s;\tilde{\delta})$. Then, $\lim_{s\rightarrow
s^*_{\tilde{\delta}}}t(s;\delta)=\infty$.
\end{lem}
\Beginproof This follows by direct substitution
of~$s^*_{\tilde{\delta}}$ in~(\ref{eq:t_s_def_2}). \Endproof
\begin{lem}\label{Lemma:Sstar_limit_}
Let $s^*:=\lim_{\tilde{\delta}\rightarrow0+}s^*_{\tilde{\delta}}$.
Then,~$s^\ast\approx-2.338$ is the first (negative) root
of~$A_i(s)$.
\end{lem}
\Beginproof By~(\ref{eq:s0_def}),~$s_0\gg1$
for~$0<\tilde{\delta}\ll1$. Therefore,
\begin{equation}\label{eq:AiryAsymptotic}
A_i(s_0)\sim
\frac1{2\sqrt{\pi}}s_0^{-1/4}e^{-\frac{2}{3}s_0^{3/2}},\qquad
B_i(s_0)\sim \frac1{\sqrt{\pi}}s_0^{-1/4}e^{\frac{2}{3}s_0^{3/2}}.
\end{equation}
Hence, by~(\ref{eq:k1_k2}),
\begin{equation}\label{eq:k1_k2_proportion}
\lim_{\tilde{\delta}\rightarrow0+}\frac{k_2}{k_1}=0.
\end{equation}
Therefore,
\[
A(s;\tilde{\delta})\sim k_1A_i(s),\qquad
s=\mathcal{O}(1),\quad0<\tilde{\delta}\ll1.
\]
Thus, the result follows. \Endproof

From Lemma~\ref{Lemma:t_s_limit_} and Lemma~\ref{Lemma:Sstar_limit_}
it follows that:
\begin{corol}
$t(s^\ast;\tilde{\delta}):=\lim_{\tilde{\delta}\rightarrow0+}t(s^\ast_{\tilde{\delta}};\tilde{\delta})=\infty$.
\end{corol}

Equation~(\ref{eq:t_s_def_2}) can be written as
\begin{equation}\label{eq:t_s_def_3}
t\left(s;\tilde{\delta}\right)=T_c^{\left(\tilde{\delta}\right)}\left[1-\frac{\frac{B_i(s)}{B_i(s_0)}}{\frac{A(s;\tilde{\delta})}{A_0}}\right],\qquad
T_c^{\left(\tilde{\delta}\right)}:=\frac{\pi}{k_1\tilde{\delta}^{1/3}}\frac{B_i(s_0)}{A_0}.
\end{equation}
By~(\ref{eq:k1_k2}) and~(\ref{eq:t_s_def_3}),
\[
T_c^{\left(\tilde{\delta}\right)}=\frac{\frac{\pi
B_i(s_0)}{k_1\tilde{\delta}^{1/3}A_0}}{\pi\left(\frac{A_0'B_i(s_0)}{(A_0)^2\tilde{\delta}^{1/3}}+A_0B_i'(s_0)\right)}.
\]
By~(\ref{eq:s0_def}) and~(\ref{eq:AiryAsymptotic}),~$B_i'(s_0)\sim
B_i(s_0)\cdot s_0^{1/2}=B_i(s_0)\cdot\delta^{-1/3}\beta_0^{1/2}$.
Therefore, since~$L=1/A$,
\[
T_c^{\left(\tilde{\delta}\right)}\sim\frac{\frac1{(A_0)^2}}{\frac{A_0'}{(A_0)^3}+\beta_0^{1/2}}=\frac{(L_0)^2}{\beta_0^{1/2}-L_0L_0'},
\]
which is the adiabatic approximation of~$T_c$, see~\cite[eq.
(3.31)]{PNLS-99}.

By~(\ref{eq:AiryAsymptotic}),~$B_i(s)/B_i(s_0)$ is exponentially
decreasing as~$s$ decreases from~$s_0$. Let us assume for simplicity
that~$A'_0\geq0$. Then,~$A(s;\tilde{\delta})>A(s_0)$.
Therefore,~$\frac{B_i(s)/B_i(s_0)}{A(s;\tilde{\delta})/A_i(s_0)}$ is
exponentially decreasing as~$s$ decreases from~$s_0$. Hence,
\[
t\left(s;\tilde{\delta}\right)\approx T_c^{(\tilde{\delta})}\approx
T_c,\qquad s^\ast<s~\mbox{and}~s_0-s\gg1.
\]
This concludes the proof of Lemma~\ref{Lemma:t_s_Lemma}. \Endproof

From Lemma~\ref{Lemma:t_s_Lemma} it follows that
\begin{equation}\label{eq:s_t_limit_Loglog}
s(t)\sim\left\{
\begin{array}{l l}
  O(s_0) & \quad 0<t<T_c,\\
  s^\ast & \quad T_c<t,\\
\end{array} \right.
\qquad as~\tilde{\delta}\rightarrow0+,
\end{equation}
see Figure~\ref{fig:t_s__s_t_DampedNLSLoglogLawWithoutNuBeta}(b). We
recall that
\[
L_t(t;\tilde{\delta})=\tilde{\delta}^{1/3}\left[k_1A_i'(s(t))+k_2B_i'(s(t))\right],
\]
see~(\ref{eq:Lt_ByAiryEq}). Hence, by~(\ref{eq:s_t_limit_Loglog}),
\[
L_t(t;\tilde{\delta})\sim
\tilde{\delta}^{1/3}\left[k_1A_i'(s^\ast)+k_2B_i'(s^\ast)\right],\qquad
T_c<t,\qquad~0<\tilde{\delta}\ll1.
\]
Therefore, by~(\ref{eq:k1_k2_proportion}),
\[
L_t(t;\tilde{\delta})\sim \tilde{\delta}^{1/3}k_1A_i'(s^\ast),\qquad
T_c<t,\quad0<\tilde{\delta}\ll1.
\]
Hence,~$L(t)$ is linear for~$T_c<t$.

By~(\ref{eq:AiryAsymptotic}) and~(\ref{eq:t_s_def_3}),
\[
k_1\tilde{\delta}^{1/3}\sim\frac{\pi}{A_0}\frac1{T_c}B_i(s_0)\sim\frac{\sqrt{\pi}}{A_0}\frac1{T_c}s_0^{-1/4}e^{\frac{2}{3}s_0^{3/2}}.
\]
{\rev{Therefore,
\begin{equation}\label{eq:LtLimitAsymtotic}
L_t(t;\tilde{\delta})\sim\frac{\sqrt{\pi}}{A_0}\frac1{T_c}s_0^{-1/4}e^{\frac{2}{3}s_0^{3/2}}A_i'(s^\ast),\qquad
T_c<t,\qquad~0<\tilde{\delta}\ll1.
\end{equation}
}}

Since~$\lim_{\tilde{\delta}\rightarrow0+}s_0=\infty$,
\[
\lim_{\tilde{\delta}\rightarrow0+}\tilde{\delta}^{1/3}k_1A'_i(s^\ast)=\infty.
\]
Therefore, the post-collapse slope of~$L(t)$ becomes infinite
as~$\tilde{\delta}\rightarrow0+$. \Endproof
\subsubsection{Simulation}
\begin{figure}[ht!]
\begin{center}
\scalebox{0.6}{\includegraphics{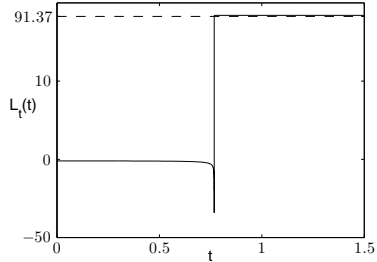}}
\caption{{\rev{Solution of the reduced
equations~(\ref{eq:DampedReducedEquation})
with~$d=1$,~$\delta=2.5\cdot10^{-4}$, and the initial
conditions~$L(0)=1$,~$L'_0=-1$, and~$\beta(0)=0.1$. Solid line
is~$L_t(t)$. Dashed line is the
prediction~(\ref{eq:LtLimitAsymtotic})}}.}
\label{fig:Lt_LoglogRedEqAsymtotic}
\end{center}
\end{figure}

{\rev{In Figure~\ref{fig:Lt_LoglogRedEqAsymtotic} we solve the
reduced equations~(\ref{eq:DampedReducedEquation}), and observe that
for~$t>T_c$,~$L_t(t)$~is indeed in excellent agreement with the
asymptotic prediction~(\ref{eq:LtLimitAsymtotic}).}}

\subsubsection{Validity of the reduced equations?}\label{subsubsec:ValidOfRE}
\begin{figure}[ht!]
\begin{center}
\scalebox{0.6}{\includegraphics{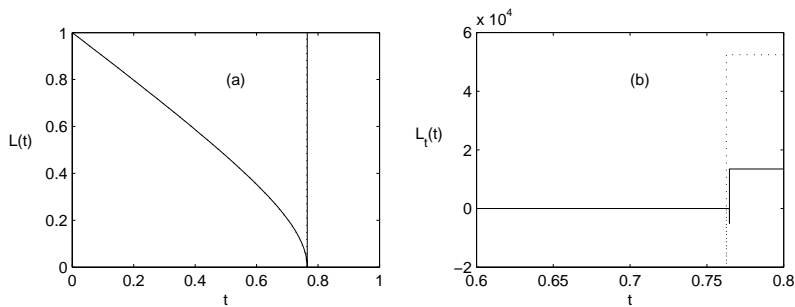}}
\caption{ Solution of the reduced
equations~(\ref{eq:beta_LeadingOrd}) [solid]
and~(\ref{eq:DampedReducedEquation}) [dots] with~$\delta=10^{-4}$,
and the initial conditions~$\beta(0)=0.1$,~$L_0=1$, and~$L_t(0)=-1$.
(a)~$L(t)$. (b)~$L_t(t)$.} \label{fig:RedEq_WithAndWithoutNuBeta}
\end{center}
\end{figure}

{\rev{In the asymptotic analysis in
Section~\ref{subsec:ReducedEqDampedNLS-loglog} we approximated the
critical NLS with the reduced equations~(\ref{eq:beta_LeadingOrd}).
Then, we approximated the reduced
equations~(\ref{eq:beta_LeadingOrd}) with the reduced
equations~(\ref{eq:DampedReducedEquation}), by
neglecting~$\nu(\beta)$. We now consider the validity of these
approximations.

In Figure~\ref{fig:RedEq_WithAndWithoutNuBeta}(a) we plot~$L(t)$ for
the solutions of the reduced equations~(\ref{eq:beta_LeadingOrd})
and~(\ref{eq:DampedReducedEquation}).
Although~$\tilde{\delta}\approx0.0019$ is not much larger
then~$\nu(\beta(0))\approx0.0016$, the two solutions are "close".
Plotting~$L_t(t)$ shows that in both cases,~$L_t(t)$ is a constant
after the collapse is arrested, see
Figure~\ref{fig:RedEq_WithAndWithoutNuBeta}(b). The "addition"
of~$\nu(\beta)$, however, decreases this constant by a factor
of~$\approx4$. Therefore, when~$\nu(\beta)$ is not neglected, the
approximation~(\ref{eq:LtLimitAsymtotic}) for~$L_t(t)$ is not
accurate, but the solution still expands linearly at a velocity that
goes to infinity as~$\delta\rightarrow0+$. }}

\begin{figure}[ht!]
\begin{center}
\scalebox{0.6}{\includegraphics{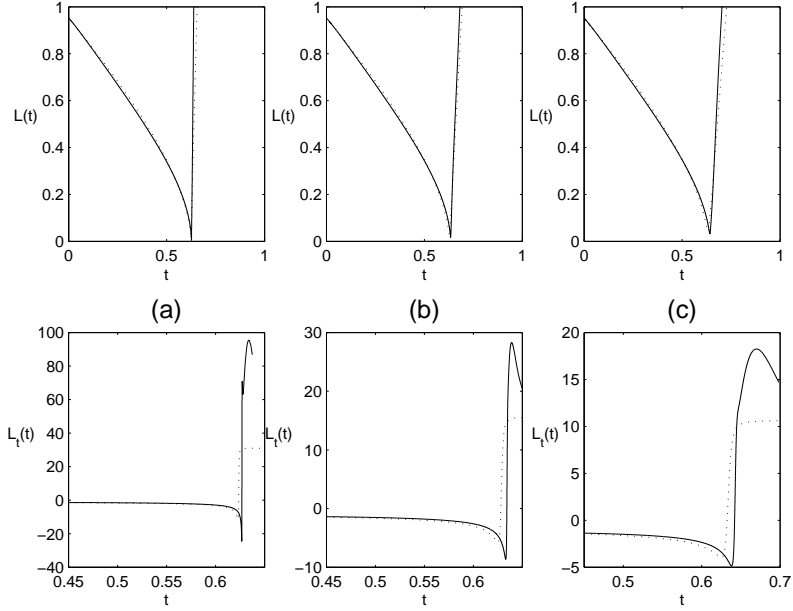}}
\caption{ Solution of the damped NLS~(\ref{eq:CriticalDampedNLS})
with~$d=1$ and the initial
condition~$\sqrt{1.05}\psi_{explicit}(t=0,r;T_c=1)$ (solid). Dashed
line is the solution of the reduced
equations~(\ref{eq:beta_LeadingOrd}) with the initial
condition~$\beta(0)\approx0.3242$,~$L(0)\approx0.9524$,
and~$L_t(0)=-1.02$.
(a)~$\delta=10^{-3}$.~(b)~$\delta=1.5\cdot10^{-3}$.~(c)~$\delta=2\cdot10^{-3}$.
Top row:~$L(t)$. Bottom row:~$L_t(t)$.}
\label{fig:DampedNLSLoglog_Sim_vs_ReducedEq_New}
\end{center}
\end{figure}

{\rev{In order to check the validity of the reduced
equations~(\ref{eq:beta_LeadingOrd}) with~$\nu(\beta)$, we solve the
damped NLS~(\ref{eq:CriticalDampedNLS}) for~$d=1$ with the initial
condition~$\sqrt{1.05}\psi_{\rm explicit}(t=0,r;T_c=1)$, for various
values of~$\delta$. Then, we extract from~$\psi$ the value of~$L(t)$
using~(\ref{eq:BeamWidthDefinition}). These NLS solutions are
compared with solutions of the reduced
equations~(\ref{eq:beta_LeadingOrd}). The initial conditions for the
reduced equations are as follow.
By~(\ref{eq:BetaDefGaussian}),~$\beta(0)\sim\frac{\|\sqrt{1.05}\psi_{explicit}\|_2^2-P_{cr}}{M}=\frac{0.05P_{cr}}{M}\approx0.3242$.
By~(\ref{eq:BeamWidthDefinition}),~$L(0)\approx0.9524$.
For~$\psi_{\rm explicit}(t=0,r;T_c=1)$,~$L_t(0)=-1$,
see~(\ref{eq:ReducedEqIC}). The multiplication by~$\sqrt{1.05}$
leads to a small change in~$L_t(0)$. We found that~$L_t(0)=-1.02$
provides the best fitting.

Figure~\ref{fig:DampedNLSLoglog_Sim_vs_ReducedEq_New} shows that
there is a good agreement between~$L(t)$ of the reduced
equations~(\ref{eq:beta_LeadingOrd}) and of the NLS. In addition, in
both cases, the post-collapse defocusing velocity increases
as~$\delta$ decreases. The curves of~$L_t(t)$ show a good agreement
when the solutions focus, but differ when the solutions defocus. In
particular,~$L_t(t)$ of the NLS solution is not a constant after the
arrest of collapse. We relate this difference to the interaction
between the expanding core and the tail, which is ignored in the
 reduced equations. This core-tail interaction did
not occur in the explicit continuation case, see
Section~\ref{subsec:ExplicitContDampedNLS}, since in that case the
power of the initial condition~$\psi_{\rm explicit}$ is equal
to~$P_{cr}$, hence there is no tail. In addition, this phenomenon
did not occur in the sub-threshold continuation (see
Section~\ref{sec:NonExplicitBeam}), since there the expansion
velocity was finite. Therefore, sufficiently close to~$T_c+$, the
expanding core had a negligible interaction with the tail.

In summary, the asymptotic analysis and numerical simulations suggest
that in the nonlinear-damping continuation of NLS solutions that undergo a loglog collapse,
the singular core of the NLS solution expands
after the singularity with a velocity that goes to infinity
as~$\delta\rightarrow0+$. The post-collapse expansion velocity is, however,
probably not linear in~$t$.}}

\section{Complex-Ginzburg Landau continuation}
\label{sec:CGL}

The two-dimensional Complex Ginzburg-Landau equation (CGL)
$$
 i \psi_t(t,x,y) + \Delta \psi + |\psi|^2 \psi
       -i \epsilon_1 \psi
        -i   \epsilon_2 \Delta \psi
        +i  \epsilon_3 |\psi|^2\psi  = 0,
$$
arises in a variety of physical problems:
 Models of  chemical turbulence,
 analysis of Poiseuille flow, Rayleigh-B\'ernard convection and Taylor-Couette flow.
Its name comes from the field of superconductivity,
 where it models  phase transitions of
materials between superconducting and non-superconducting phases.

In~\cite{CGL-99}, Fibich and Levy showed that as
$\epsilon_1,\epsilon_2,\epsilon_3 \longrightarrow 0$, the collapse
dynamics is governed, to leading order, by the reduced
equations~\eqref{eq:DampedReducedEquation} with
$$
\tilde\delta = \frac{2 \Pcr}{M} (\epsilon_2+ 2 \epsilon_3).
$$
Therefore, Continuation Results~\ref{Proposition:DampedNLSWeakSol}
and~\ref{conject:DampedNLS_Conject} hold also for the CGL
continuation of the critical NLS.

\section{Continuation in the linear \schro~equation}
\label{sec:SingularityInLS}

 It is well known that in the linear
\schro~equation, under the geometrical-optics approximation, a
focused input beam becomes singular at the focal point. When,
however, diffraction is not neglected, the focused beam does not
collapse to a point, but rather narrows down to a positive
\textit{diffraction-limited width}, and then spreads out with
further propagation. Therefore, diffraction can play the role of
"viscosity" in the continuation of singular geometrical-optics
linear solutions. In what follows, we compare this continuation with
those in the nonlinear case.

Consider the d-dimensional linear \schro~equation
\begin{subequations}\label{eq:DimensionlessScrodingerEquation}
\begin{equation}
2ik_0\psi_t(\X,t)+\Delta\psi=0,\qquad \X=(x_1,...,x_d),
\end{equation}
with a focused Gaussian initial condition
\begin{equation}
\psi_0(\X)=e^{-r^2/2}e^{-ik_0r^2/2F},\qquad r=|\X|,
\end{equation}
\end{subequations}
where~$F>0$ is the focal point. We can look for a solution of the
form
\begin{equation}\label{eq:GenericInitialCondition}
\psi(\X,t)=A(\X,t)e^{ik_0S(\X,t)},
\end{equation}
where A and~$S$ are real. Substitution in
equation~(\ref{eq:DimensionlessScrodingerEquation}) gives
\begin{subequations}\label{eq:SeperatedAmpPhaseSolution}
\begin{equation}\label{eq:SeperatedAmpPhaseSolution_real}
\left(\nabla S\right)^2+2S_t-\frac1{k_0^2}\frac{\Delta A}{A}=0,
\end{equation}
and
\begin{equation}\label{eq:SeperatedAmpPhaseSolution_imag}
\left(A^2\right)_t+\nabla S\cdot\nabla\left(A^2\right)+A^2\Delta
S=0,
\end{equation}
\end{subequations}
with initial conditions
\begin{equation}\label{eq:SeperatedIC}
S(\X,0)=-\frac{r^2}{2F},\qquad A^2(\X,0)=e^{-r^2}.
\end{equation}
Since~$k_0\gg 1$, we can apply the geometrical-optics approximation,
and neglect the diffraction term~$\Delta A$. In this case,
equation~(\ref{eq:SeperatedAmpPhaseSolution_real}) becomes
\begin{equation}\label{eq:EikonalEquation}
\left(\nabla S\right)^2+2S_t=0,
\end{equation}
while equation~(\ref{eq:SeperatedAmpPhaseSolution_imag}) remains
unchanged.

The solution of equations~(\ref{eq:SeperatedAmpPhaseSolution_imag},
~\ref{eq:EikonalEquation}), subject to the initial
conditions~(\ref{eq:SeperatedIC}), is given by
\[
S=\frac{r^2}{2}\frac{L_t}{L},\qquad
A^2(\X,t)=\frac1{L^d(t)}e^{-\frac{r^2}{L^2(t)}},\qquad
L(t)=1-\frac{t}{F}.
\]
Therefore, under the geometrical-optics approximation, the solution
of~(\ref{eq:DimensionlessScrodingerEquation}) is given by
\begin{subequations}\label{eq:GaussianLinearScrodingerGO}
\begin{equation}\label{eq:GaussianLinearScrodingerGO_profile}
\psi^{go}_{lin}(t,\X)=\frac1{L_{go}^{d/2}(t)}e^{-\frac1{2}\frac{r^2}{L_{go}^2(t)}}e^{ik_0\frac{r^2}{2}\frac{L'_{go}(t)}{L_{go(t)}}}\quad,\qquad
0\leq t< F,
\end{equation}
where
\begin{equation}\label{eq:GaussianLinearScrodingerGO_width}
L_{go}(t)=1-\frac{t}{F}.
\end{equation}
\end{subequations}
Since~$\lim_{t\rightarrow
F}|\psi^{go}_{lin}|^2=\|\psi_0\|^2_2\cdot\delta(\X)$, the
geometrical-optics solution becomes singular at the focal
point~$t=F$.

Equation~(\ref{eq:DimensionlessScrodingerEquation}) can also be
solved exactly, (i.e., without making the geometrical-optics
approximation), yielding
\begin{equation}\label{eq:GaussianBeamExplicitSolution_psi}
\psi_{lin}(t,\X)=\frac1{L^{d/2}(t)}e^{-\frac1{2}\frac{r^2}{L^2(t)}}e^{ik_0\frac{r^2}{2}\frac{L_t}{L}+i\tau(t;k_0)},
\end{equation}
where
\begin{equation}\label{eq:GaussianBeamExplicitSolution_width}
L(t;k_0)=\sqrt{\frac1{F}\frac{(t-t_{min})^2}{t_{min}}+L^2_{min}},\qquad
\tau(t;k_0)=-\frac{d}{2}\left[\mbox{atan}\left(\frac{t-t_{min}}{L_{min}\sqrt{F\cdot
t_{min}}}\right)+\mbox{atan}\left(\frac{k_0}{F}\right)\right],
\end{equation}
and
\[
t_{min}=\frac{F}{1+F^2/k_0^2},\qquad
L_{min}=\frac{F}{\sqrt{F^2+k_0^2}}.
\]
Since~$L(t;k_0)$ does not shrink to zero at any~$t>0$,~$\psi_{lin}$
exists for all~$0\leq t<\infty$, and in particular for~$t>F$.

It is easy to verify that the limiting width of~$\psi_{lin}$ is
given by
\begin{equation}\label{eq:LimitingK0BeamWidth}
\lim_{k_0\rightarrow\infty}L(t;k_0)=|L_{go}(t)|.
\end{equation}
In addition, since~$\lim_{k_0\rightarrow\infty}t_{min}=F$
and~$\lim_{k_0\rightarrow\infty}L_{min}=0$, then
\begin{equation}\label{eq:LinearSchro_TauLimit}
\lim_{k_0\rightarrow\infty}\tau(t;k_0)=\left\{
\begin{array}{l l}
  0 & \quad 0\leq t<F,\\
  -\frac{d}{2}\cdot\pi & \quad F<t<\infty.
\end{array} \right.
\end{equation}
Therefore,
\begin{lem}
\begin{equation}\label{eq:GeometricalOpticsLimitSolution}
\lim_{k_0\rightarrow\infty}\psi_{lin}(t,\X)=\left\{
\begin{array}{l l}
  \psi_{lin}^{go}(t,\X) & \quad 0\leq t<F,\\
  \left(\psi_{lin}^{go}\right)^\ast(2F-t,\X)e^{-i\frac{d}{2}\cdot\pi} & \quad F<t<\infty.
\end{array} \right.
\end{equation}
\end{lem}
Hence,~$\lim_{k_0\rightarrow\infty}\psi$ coincides
with~$\psi^{go}_{lin}$ before the focal point. Beyond the focal
point, there is a bounded jump in the limiting phase, and the
continuation is symmetric with respect to~$T_c=F$. This symmetry is
to be expected, since the linear continuation is invariant
under~(\ref{eq:NLSInvariantTransformations}).

By~(\ref{eq:LinearSchro_TauLimit}), the limiting phase is unique,
both before and after the singularity. This is the opposite from the
nonlinear case, where the limiting phase is non-unique beyond the
singularity. We thus conclude that
\begin{conc}
The post-collapse non-uniqueness of the phase is a nonlinear
phenomena.
\end{conc}
\subsection{Simulations}
\begin{figure}[ht!]
\begin{center}
\scalebox{0.6}{\includegraphics{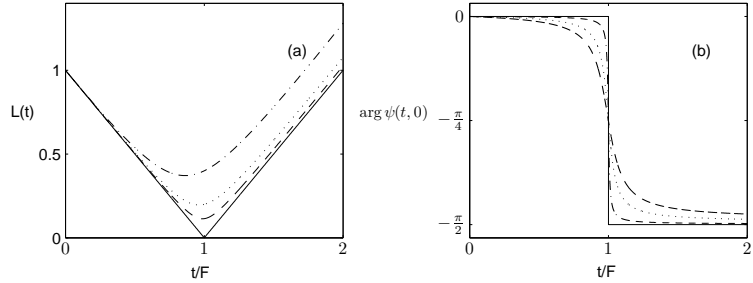}}
\caption{ Solution of the linear
\schro~(\ref{eq:DimensionlessScrodingerEquation}) with~$d=1$
and~$F=4$. (a): Width~$L(t)$,
see~(\ref{eq:GaussianBeamExplicitSolution_width}), for~$k_0=10$
(dashes-dots),~$k_0=20$ (dots) and~$k_0=35$ (dashes). Solid line
is~$|L_{go}(t)|=\left|1-\frac{t}{F}\right|$.
(b):~$\arg\psi(t,0;k_0)$ as a function of~$t$ for~$k_0=50$
(dashes),~$k_0=100$ (dots) and~$k_0=500$ (dashes-dots). Solid line
is~(\ref{eq:LinearSchro_TauLimit}).}\label{fig:BeamWidthLinearSchrodinger_GO}
\end{center}
\end{figure}

In order to illustrate these results numerically, in
Figure~\ref{fig:BeamWidthLinearSchrodinger_GO}(a) we
plot~$L(t;k_0)$, and observe that it approaches~$|L_{go}(t)|$
as~$k_0\rightarrow\infty$, both before and after the singularity
point at~$t=F$. In Figure~\ref{fig:BeamWidthLinearSchrodinger_GO}(b)
we plot~$\tau(t;k_0)=\arg\psi(t,0;k_0)$ as a function of~$t$, and
observe that as~$k_0$ increases,~$\tau(t;k_0)$ approaches the step
function~(\ref{eq:LinearSchro_TauLimit}).

\section{Discussion}
   \label{sec:FinalRemarks}

In this study, we presented several "old" continuations and four
novel continuations of NLS solutions. At this stage, it is not clear
whether any of these continuations is the ``physical" one. In fact,
it is possible that there is no "universal continuation", i.e, that
different physical setups calls for different continuations.

This study suggests that all continuations share the
property that the post-collapse phase becomes non-unique. Indeed,
this property follows from the fact that the phase of the singular
solution becomes infinite at the singularity,  and is thus
independent of the specific continuation which is used. Therefore,
even without knowing the "correct" continuation, we can conclude
that interactions between post-collapse filaments are chaotic.

Since the Bourgain-Wang blowup solutions are unstable, they are
typically classified as ``non-generic'' solutions. This study shows,
however, that these solutions are ``generic'', in the sense that
they arise as the sub-threshold power continuation of NLS solutions,
both before and after the singularity.

It is instructive to compare this study with~\cite{Merle-Raphael-Szeftel}. In~\cite{Merle-Raphael-Szeftel},
Merle, Raphael and Szeftel showed that the Bourgain-Wang blowup solutions lie on the boundary of an $H^1$ open set of global solutions that scatter forward and backward in time, and also on the boundary of an $H^1$ open set of solutions that undergo a loglog blowup in finite time. This result follows immediately from the
reduced equations~(\ref{eq:ReducedEqGaussian})--(\ref{eq:ReducedEqCont}) of the sub-threshold power continuation, see Section~\ref{subsec:ProofOfProp1}, since
\begin{itemize}
  \item The Bourgain-Wang blowup solutions correspond to $\beta(0) = 0$.
  \item  For any~$\beta(0)<0$, $L(t)$ remains strictly positive and satisfies $\lim_{t \to \pm \infty} L(t) = 0$.
  \item For any $\beta(0)>0$, $L(t)$ goes to zero at some finite~$T_c$, at the loglog rate.
\end{itemize}
An interesting difference between the approaches of this study and~\cite{Merle-Raphael-Szeftel}, is that Merle, Raphael and Szeftel start from the Bourgain-Wang solution at the singularity time~$T_c$, and then find a smooth deformation such that the deformed solutions belong to the above two open sets on either side of the Bourgain-Wang solution.
In this study, we start from a generic initial condition~$\psi_0 = K \cdot F({\bf x})$, and
obtain the Bourgain-Wang solution as $K \longrightarrow K_{th}$.\footnote{Of course,
another difference is that, unlike this study, the results of~\cite{Merle-Raphael-Szeftel} are rigorous.}

The `global'' picture that emerges from this study is as follows. Consider a
stable singular solution of the critical NLS that undergoes a loglog collapse.
Since the singular core~$\psi_{R^{(0)}}$ approaches a~$\delta$-function with
power~$\Pcr$, it can be continued with a $\delta$-function filament with power~$\Pcr$.
Such a filament singularity can occur when the collapse-arresting mechanism
leads to focusing-defocusing oscillations. Since this is the generic effect of conservative perturbations of the critical NLS (such as nonlinear saturation), see~\cite[Section 4.1.2]{PNLS-99}, we expect that
{\em continuations that are based on conservative perturbations of the NLS will lead to a filament singularity}.

When the collapse-arresting mechanism is non-conservative (e.g.,
nonlinear damping), the solution defocuses (scatters) after its
collapse has been arrested, since its power gets below~$\Pcr$.
Therefore, we expect that {\em continuations that are based on
non-conservative perturbation of the NLS will lead to a point
singularity.} In addition, the same arguments as in the proof of
Continuation Result~\ref{Proposition:WeakSolutionSymmetry} suggest
that {\em non-conservative continuations of solutions that undergo a
loglog collapse have an infinite-velocity expanding{\rev{core}}}.

In Continuation Result~\ref{Proposition:WeakSolutionSymmetry} we saw that when the continuation leads to a point
singularity and is time-reversible, the continuation is symmetric with respect to~$T_c$ (Property 1).
For this to occur, however, the continuation should be conservative (in order to be time reversible),
yet it should not lead to focusing-defocusing oscillations. While this
holds for the sub-threshold power and the shrinking-hole continuations, it is not
expected to hold for conservative perturbations of the NLS, which generically lead to a focusing-defocusing oscillations, hence to a filament singularity. Therefore, we expect that {\em continuations which are based on perturbations of the NLS equation are asymmetric with respect to~$T_c$.} Hence, we believe that Property~1 is non-generic.

\subsection*{Acknowledgment}
We acknowledge useful discussions with Frank Merle. This research
was partially supported by grant~$1023/08$ from the Israel Science
Foundation (ISF).


\begin{thebibliography}{10}

\bibitem{Kelley-65}
P.L. Kelley.
\newblock Self-focusing of optical beams.
\newblock {\em Phys. Rev. Lett.}, 15:1005--1008, 1965.

\bibitem{Merle-92a}
F.~Merle.
\newblock On uniqueness and continuation properties after blow-up time of
  self-similar solutions of nonlinear {S}chr\"odinger equation with critical
  exponent and critical mass.
\newblock {\em Comm. Pure Appl. Math.}, 45:203--254, 1992.

\bibitem{Merle-92b}
F.~Merle.
\newblock Limit behavior of saturated approximations of nonlinear
  {S}chr\"{o}dinger equation.
\newblock {\em Comm. Math. Phys.}, 149:377--414, 1992.

\bibitem{Merle-Raphael-Szeftel}
{F.} {M}erle, {P.} {R}aphael, and {J.} {S}zeftel.
\newblock The instability of {B}ourgain-{W}ang solutions for the {L}$^2$
  critical {N}{L}{S}.
\newblock {\em preprint}.

\bibitem{Tao-2009}
T.~Tao.
\newblock Global existence and uniqueness results for weak solutions of the
  focusing mass-critical nonlinear {S}chr\"odinger equation.
\newblock {\em Analysis and PDE}, 2:61--81, 2009.

\bibitem{Stinis-2010}
P.~Stinis.
\newblock Numerical computation of solutions of the critical nonlinear
  {S}chr\"odinger equation after the singularity.
\newblock {\em preprint}, 2010.

\bibitem{Bourgain-97}
J.~Bourgain and W.~Wang.
\newblock Construction of blowup solutions for the nonlinear {S}chr\"odinger
  equation with critical nonlinearity.
\newblock {\em Ann. Scuola Norm. Sup. Pisa Cl. Sci.}, 25:197--215, 1997.

\bibitem{Malkin-1993}
V.~Malkin.
\newblock On the analytical theory for stationary self-focusing of radiation.
\newblock {\em Physica D}, 64:251--266, 1993.

\bibitem{PNLS-99}
G.~Fibich and G.C. Papanicolaou.
\newblock Self-focusing in the perturbed and unperturbed nonlinear
  {S}chr\"odinger equation in critical dimension.
\newblock {\em Siam Appl. Math.}, 60:183--240, 1999.

\bibitem{Sulem-99}
C.~Sulem and P.L. Sulem.
\newblock {\em The Nonlinear {S}chr\"{o}dinger Equation}.
\newblock Springer, New-York, 1999.

\bibitem{Strauss-89}
W.~Strauss.
\newblock {\em Nonlinear wave equation}.
\newblock American Mathematical Society, Providence, R.I, 1989.

\bibitem{Weinstein-83}
M.I. Weinstein.
\newblock Nonlinear {S}chr\"{o}dinger equations and sharp interpolation
  estimates.
\newblock {\em Comm. Math. Phys.}, 87:567--576, 1983.

\bibitem{Merle-Raphael-2003}
F.~Merle and P.~Raphael.
\newblock Sharp upper bound on the blowup rate for the critical nonlinear
  {S}chr\"{o}dinger equation.
\newblock {\em Geom. Funct. Anal}, 13:591--642, 2003.

\bibitem{Merle-04}
F.~Merle and P.~Raphael.
\newblock On universality of blow-up profile for {$L\sp 2$} critical nonlinear
  {S}chr\"odinger equation.
\newblock {\em Invent. Math.}, 156:565--672, 2004.

\bibitem{Merle-05}
F.~Merle and P.~Raphael.
\newblock Blow-up dynamics and upper bound on the blow-up rate for the critical
  nonlinear {S}chr\"odinger equation.
\newblock {\em Ann. of Math.}, 161:157--222, 2005.

\bibitem{Merle-05b}
F.~Merle and P.~Raphael.
\newblock Profiles and quantization of the blow-up mass for critical nonlinear
  {S}chr\"odinger equation.
\newblock {\em Commun. Math. Phys.}, 253:675--704, 2005.

\bibitem{Merle-06}
F.~Merle and P.~Raphael.
\newblock On a sharp lower bound on the blow-up rate for the {$L\sp 2$}
  critical nonlinear {S}chr\"odinger equation.
\newblock {\em J. Amer. Math. Soc.}, 19:37--90, 2006.

\bibitem{Merle-06b}
F.~Merle and P.~Raphael.
\newblock On one blow up point solutions to the critical nonlinear
  {S}chr\"odinger equation.
\newblock {\em J. Hyperbolic Differ. Eq.}, 2:919--962, 2006.

\bibitem{Raphael-05}
P.~Raphael.
\newblock Stability of the log-log bound for blow up solutions to the critical
  non linear {S}chr\"odinger equation.
\newblock {\em Math. Ann.}, 331:577--609, 2005.

\bibitem{Papanicolaou-88_1}
{B.} {L}eMesurier, {C.} {S}ulem, {G.} {P}apanicolaou, and {P.L.} {S}ulem.
\newblock Local structure of the self-focusing singularity of the nonlinear
  {S}chr\"{o}dinger equation.
\newblock {\em Phys. D}, 38:210--226, 1988.

\bibitem{Papanicolaou-88_2}
{C.} {S}ulem, {G.} {P}apanicolaou, {M.} {L}andman, and {P.L.} {S}ulem.
\newblock Rate of blowup for solutions of the nonlinear {S}chr\"{o}dinger
  equation at critical dimension.
\newblock {\em Phys. Rev. A}, 32:3837--3843, 1988.

\bibitem{Fraiman-1985}
G.~Fraiman.
\newblock Asymptotic stability of manifold of self-similar solutions in
  self-focusing.
\newblock {\em Phys. JETP}, 61:228--233, 1985.

\bibitem{Cazenave-90}
T.~Cazenave and F.B. Weissler.
\newblock The {C}auchy problem for the critical nonlinear {S}chr\"{o}dinger
  equation in {$H^s$}.
\newblock {\em Nonlinear Anal.}, 14:807--836, 1990.

\bibitem{Gprofile-05}
{G.} {F}ibich, {N.} {G}avish, and {X.P.} {W}ang.
\newblock New singular solutions of the nonlinear {S}chr\"odinger equation.
\newblock {\em Physica D}, 211:193--220, 2005.

\bibitem{SC_rings-07}
{G.} {F}ibich, {N.} {G}avish, and {X.P.} {W}ang.
\newblock Singular ring solutions of critical and supercritical nonlinear
  {S}chr\"{o}dinger equations.
\newblock {\em Physica D}, 231:55--86, 2007.

\bibitem{Baruch_Fibich_Gavish:2009}
{G.} {B}aruch, {G.} {F}ibich, and {N.} {G}avish.
\newblock Singular standing ring solutions of nonlinear partial differential
  equations.
\newblock {\em Physica D}, 2009.

\bibitem{Markowich-04}
{W.} {B}ao, {D.} {J}aksch, and {P.A.} {M}arkowich.
\newblock Three-dimension simulation of jet formation in collapsing
  condensates.
\newblock {\em J. Phys. B: At. Mol. Opt. Phys.}, 37:329--343, 2004.

\bibitem{Fibich-2001}
G.~Fibich.
\newblock Self-focusing in the damped nonlinear {S}chr\"odinger equation.
\newblock {\em SIAM Appl. Math.}, 61:1680--1705, 2001.

\bibitem{Sulem-Passot-2005}
{T.} {P}assot, {C.} {S}ulem, and {P.L.} {S}ulem.
\newblock Linear versus nonlinear dissipation for critical {N}{L}{S} equation.
\newblock {\em Physica D}, 203:167--184, 2005.

\bibitem{Antoneli-Sparber-2010}
P.~Antonelli and C.~Sparber.
\newblock Global well-posedness for cubic {N}{L}{S} with nonlinear damping.
\newblock {\em Comm. PDE}, 35:4832--4845, 2010.

\bibitem{Merle-private-11}
F.~Merle.
\newblock Private communication.
\newblock 2011.

\bibitem{Albright-Gavathas-86}
J.~R. Albright and E.~P. Gavathas.
\newblock Integrals involving {A}iry functions.
\newblock {\em Phys. A: Math. Gen}, 19:2663--2665, 1986.

\bibitem{CGL-99}
G.~Fibich and D.~Levy.
\newblock Self-focusing in the complex {G}inzburg-{L}andau limit of the
  critical nonlinear {S}chr\"odinger equation.
\newblock {\em Phys. Lett. A}, 249:286--294, 1998.

\end{thebibliography}

\end{document}